%% file: Arxiv_Published.tex
\documentclass[12pt,a4paper]{article}
\usepackage{textcomp}
\usepackage{listings}
\usepackage[space=true]{accsupp}
\newcommand{\copyablespace}{\BeginAccSupp{method=hex,unicode,ActualText=00A0}\EndAccSupp{}}
\usepackage[T1]{fontenc}
\usepackage[english]{babel}
\usepackage{hyperref}
\usepackage{graphicx}
\usepackage{amsmath}
\usepackage{amssymb}
\usepackage{mathtools}
\usepackage[table,dvipsnames]{xcolor}
\usepackage{fancyhdr}
\usepackage{longtable}
\usepackage{multicol}
\usepackage{enumerate}
\usepackage{enumitem}
\setlist[itemize]{leftmargin=5.5mm}
\usepackage{multirow}
\usepackage{natbib}
\usepackage{caption}
\usepackage{subcaption}
\usepackage{xurl}
\usepackage[utf8]{inputenc}
\usepackage{tikz}
\usepackage{pgfplots}
\pgfplotsset{compat=newest}

\pgfplotsset{
    box plot/.style={
        /pgfplots/.cd,
        black,
        only marks,
        mark=-,
        mark size=\pgfkeysvalueof{/pgfplots/box plot width},
        /pgfplots/error bars/y dir=plus,
        /pgfplots/error bars/y explicit,
        /pgfplots/table/x index=\pgfkeysvalueof{/pgfplots/box plot x index},
    },
    box plot box/.style={
        /pgfplots/error bars/draw error bar/.code 2 args={%
            \draw  ##1 -- ++(\pgfkeysvalueof{/pgfplots/box plot width},0pt) |- ##2 -- ++(-\pgfkeysvalueof{/pgfplots/box plot width},0pt) |- ##1 -- cycle;
        },
        /pgfplots/table/.cd,
        y index=\pgfkeysvalueof{/pgfplots/box plot box top index},
        y error expr={
            \thisrowno{\pgfkeysvalueof{/pgfplots/box plot box bottom index}}
            - \thisrowno{\pgfkeysvalueof{/pgfplots/box plot box top index}}
        },
        /pgfplots/box plot
    },
    box plot top whisker/.style={
        /pgfplots/error bars/draw error bar/.code 2 args={%
            \pgfkeysgetvalue{/pgfplots/error bars/error mark}%
            {\pgfplotserrorbarsmark}%
            \pgfkeysgetvalue{/pgfplots/error bars/error mark options}%
            {\pgfplotserrorbarsmarkopts}%
            \path ##1 -- ##2;
        },
        /pgfplots/table/.cd,
        y index=\pgfkeysvalueof{/pgfplots/box plot whisker top index},
        y error expr={
            \thisrowno{\pgfkeysvalueof{/pgfplots/box plot box top index}}
            - \thisrowno{\pgfkeysvalueof{/pgfplots/box plot whisker top index}}
        },
        /pgfplots/box plot
    },
    box plot bottom whisker/.style={
        /pgfplots/error bars/draw error bar/.code 2 args={%
            \pgfkeysgetvalue{/pgfplots/error bars/error mark}%
            {\pgfplotserrorbarsmark}%
            \pgfkeysgetvalue{/pgfplots/error bars/error mark options}%
            {\pgfplotserrorbarsmarkopts}%
            \path ##1 -- ##2;
        },
        /pgfplots/table/.cd,
        y index=\pgfkeysvalueof{/pgfplots/box plot whisker bottom index},
        y error expr={
            \thisrowno{\pgfkeysvalueof{/pgfplots/box plot box bottom index}}
            - \thisrowno{\pgfkeysvalueof{/pgfplots/box plot whisker bottom index}}
        },
        /pgfplots/box plot
    },
    box plot median/.style={
        /pgfplots/box plot,
        /pgfplots/table/y index=\pgfkeysvalueof{/pgfplots/box plot median index}
    },
    box plot width/.initial=1em,
    box plot x index/.initial=0,
    box plot median index/.initial=1,
    box plot box top index/.initial=2,
    box plot box bottom index/.initial=3,
    box plot whisker top index/.initial=4,
    box plot whisker bottom index/.initial=5,
}

\newcommand{\boxplot}[2][]{
    \addplot [box plot median,#1] table {#2};
    \addplot [forget plot, box plot box,#1] table {#2};
    \addplot [forget plot, box plot top whisker,#1] table {#2};
    \addplot [forget plot, box plot bottom whisker,#1] table {#2};
}

\usepackage{framed}
\usepackage{booktabs}
\usepackage{caption}
\usepackage{float}
\usepackage{titlesec}
\usepackage{capt-of}
\definecolor{darkgreen}{rgb}{0.1, 0.5, 0.2}
\usepackage{array}
\usepackage{arydshln}
\newtheorem{definition}{Definition}

\newtheorem{proposition}{Proposition}
\newtheorem{remark}{Remark}

\definecolor{gray2}{rgb}{0.6,0.6,0.6}
\definecolor{lightgray2}{rgb}{0.8,0.8,0.8}
\title{Optimal sequences for pairwise comparisons: the graph of graphs approach}
\author{Zsombor Sz\'adoczki$^{1,2,*}$, S\'andor Boz\'oki$^{1,2}$}
\date{}
\begin{document}
\pagenumbering{arabic}

\maketitle
\begin{center}
$*$ Corresponding author, 1111 Kende u. 13-17., Budapest, Hungary;\\ Email: szadoczki.zsombor@sztaki.hu\\
\bigskip
$^{1}$ Research Group of Operations Research and Decision Systems, \\
Research Laboratory on Engineering \& Management Intelligence \\
HUN-REN Institute for Computer Science and Control (HUN-REN SZTAKI), 1111 Kende u. 13-17., Budapest, Hungary;\\ Email: szadoczki.zsombor@sztaki.hu, bozoki.sandor@sztaki.hu\\
\bigskip
$^{2}$ Department of Operations Research and Actuarial Sciences \\
Corvinus University of Budapest, 1093 Fővám tér 8., Budapest, Hungary \\

\end{center}



\begin{abstract}

\noindent
In preference modelling, it is essential to determine the number of questions and their arrangements to ask from the decision maker. We focus on incomplete pairwise comparison matrices, and provide the optimal filling in patterns, which result in the closest (LLSM) weight vectors on average to the complete case for at most six alternatives and for all possible number of comparisons, when the underlying representing graph is connected. These results are obtained by extensive numerical simulations with large sample sizes. Many optimal filling structures resulted in optimal filling in sequences---one optimal case can be reached by adding a comparison to a previous one---which are presented on graph \color{black} of graphs. The star graph is revealed to be optimal among spanning trees, while the optimal graphs are always close to bipartite ones. Regular graphs also correspond to optimal cases, furthermore regularity is important for all optimal graphs, as the degrees of different vertices are always as close to each other as possible. Besides applying optimal filling structures in given decision making problems, practitioners can utilize the optimal filling sequences in the cases, when the decision maker can abandon the problem at any period of the process (e.g., in online questionnaires).

\end{abstract}


\renewcommand{\baselinestretch}{1.24} \normalsize

\newpage

\section{Introduction}
\label{sec:1}

The concept of pairwise comparisons \citep{Thurstone1927} is fundamental both in preference modelling and Multicriteria Decision Making (MCDM) \citep{Triantaphyllou2000}. These comparisons are frequently placed into so-called pairwise comparison matrices (PCMs), which are the basis of the Analytic Hierarchy Process (AHP) \citep{Saaty1977,Saaty1980}. Incompleteness (the absence of some comparisons) occurs quite often in practical problems \citep{BozokiCsatoTemesi,Temesi2023}, as well as in theoretical questions \citep{Fedrizzi2007,Bozoki2010,CsatoRonyai2016,Kulakowski2020}. In connection with decision making problems, one major source of missing data is the lack of willingness or time of the decision maker, as completing all comparisons---especially in the case of many different levels, criteria, and alternatives---can be exhausting and lingering \citep{FedrizziGiove2013,Bozoki2020}.

We would like to underline that the aim of our research is not to encourage decision makers to make less comparisons or decision analysts to ask fewer questions, although it is one of the goals of several current studies \citep{Rezaei2016,Abastante2019,Duleba2022} to unburden the decision makers by reducing the number of comparisons in practical decision making. \color{black} However, we would like to provide the sequence of questions for the analysts, which ensures that whenever the decision maker stops answering the questions, the calculated preferences are in some sense the closest to the decision makers' real preferences.

The arrangement of comparisons, which has a crucial effect on the results, is often represented by graphs \citep{Gass1998}. In this paper, we are the first to provide the optimal filling in patterns of incomplete pairwise comparison matrices, which on average produce the (both cardinally and ordinally) closest weight vectors to the complete case, for at most six alternatives (criteria) ($n$) for  each \color{black} possible number of comparisons ($e$), when the respective graph is connected. These optimal patterns for the examined $(n,e)$ pairs are significant findings of this paper themselves, however, they result in (partial) optimal filling in sequences, which can be instrumental in the case of such problems (e.g., online questionnaires), where the decision makers can abandon the problem at any period of the process to always be as close to the decision makers' preferences as possible.

These kind of problems are often present in the case of large-scale group decision making \citep{Duleba2012,Tang2021,LI2022,Liangetal2023}, or when several different experts' comparisons should be evaluated from different fields as well \citep{FRANCISOLIVIERO2021}, however, the results of this paper are not limited to group decision making problems. \color{black}

In the analysis of filling in sequences, the focus of the paper, but also in structural analysis of graphs and graph sequences in general, graph \color{black}  of graphs is a convenient and efficient tool for research and visualization, too. 
 Nodes \color{black} of a graph of graphs \color{black} are graphs, and there is an edge \color{black} between two nodes \color{black} (=graphs) if the associated graphs are in a specified relation, e.g., they can be drawn from each other by adding or deleting an edge. 
Depending on the specification of the relation, several graphs \color{black} of graphs have been investigated, see for instance \cite{Lovasz1977}. Another remarkable graph \color{black} of graphs is the Petersen family of seven graphs, including the Petersen graph itself \citep{HashimotoNikkuni2013}. The graph \color{black} of graphs by \cite{Mesbahi2002} is motivated by the evolution of graphs in a dynamic system.


It is worth noting that the term ‘neighbouring graphs' in \cite{Lovasz1977} is used synonymously for ‘there is an edge \color{black} between two graphs'. Analogously, ‘reachable' in  \cite{Mesbahi2002} means that there is a path \color{black} between two graphs. We use the concept of graph \color{black} of graphs to visualize our findings throughout the paper. To make it easier to follow, the graphs of graphs and their components, from now on, are referred to using capital letters in the study (e.g., GRAPH, EDGE, NODE, PATH, etc.) distinguishing them from the graphs (NODEs) themselves. \color{black}
  
The rest of the paper is organized as follows. Section \ref{sec:2} presents the fundamental concepts and definitions regarding PCMs and their graph representation. The methodology of the applied simulations and the related probability theoretical reasoning are detailed in Section \ref{sec:3}, while Section \ref{sec:4} contains the results, the optimal filling in sequences for the examined cases. Finally, Section \ref{sec:5} concludes and raises research questions for the future.

\section{Basic concepts: PCMs and their graph representation}
\label{sec:2}

Pairwise comparisons are the core of ranking, sports competitions, as well as many statistics and decision making techniques \citep{DavidsonFarquhar1976,Csato2021}. We focus on pairwise comparison matrices (PCMs) which are used in the Analytic Hierarchy Process (AHP) MCDM methodology to evaluate alternatives according to a criterion, as well as to determine the importance of the different criteria. However, our results can be beneficial in a wider range.

\begin{definition}[Pairwise comparison matrix (PCM)]
Let us denote the number of criteria (alternatives) in a decision problem by $n$. The $n\times n$ matrix $A=[a_{ij}]$ is called a pairwise comparison matrix, if it is positive ($a_{ij}>0$ for all  $ i $  and  $ j$) and reciprocal ($1/a_{ij}  = a_{ji}$ for all $ i $ and $ j$).
\end{definition}

The element $a_{ij}$ of a PCM shows how many times item $i$ is better/stronger/more important than item $j$. However, when a decision maker fills in all $n(n-1)/2$ elements (the elements above the principal diagonal, because of the reciprocity) there can be some kind of contradiction, a certain inconsistency in the PCM.

\begin{definition}[Consistent PCM]
A PCM is said to be consistent if  $a_{ik}=a_{ij}a_{jk} \hspace{0.2cm} \forall i,j,k$. If a PCM is not consistent, then it is called inconsistent.
\end{definition}

Naturally, there are several degrees of inconsistency, which leads to the deeply analyzed problem of different inconsistency indices \citep{Brunelli2018}, their properties \citep{Brunelli2017}, and the appropriate recommended thresholds \citep{Amenta2020,Agoston2022}. Although, many measures have been proposed, the most widely used one is probably still Saaty's Consistency Ratio (CR) \citep{Saaty1977}.

\begin{definition}[Consistency Ratio (CR)]
The CR of an $n\times n$ PCM $A$ is defined as follows:
\begin{equation}
\label{eq:1}
    CR=\frac{CI}{RI},
\end{equation}
where CI stands for Consistency Index, that is:
\begin{equation}
\label{eq:2}
    CI=\frac{\lambda_{\max}-n}{n-1},
\end{equation}
where $\lambda_{\max}$ is the principal eigenvalue of the matrix $A$, and RI is the Random Index, which is the average CI obtained from a sufficiently large set of randomly generated PCMs of size $n$.
\end{definition}

Probably the two most commonly used techniques to calculate a weight vector (prioritization vector) from a PCM that shows the importance of compared items, are the logarithmic least squares (LLSM) \citep{Crawford1985} and the eigenvector (EV) \citep{Saaty1977} methods.

\begin{definition}[Logarithmic Least Squares Method (LLSM)]
Let $A$ be an $n\times n$ PCM. The weight vector $w$ of $A$ determined by the LLSM is given as follows:
\begin{equation}
\label{eq:3}
\min_{w}    \sum_{i=1}^n\sum_{j=1}^n \left(\ln(a_{ij})-\ln\left(\frac{w_i}{w_j}\right)\right)^2 ,
\end{equation}
where $w_i$ is the $i$th coordinate of $w$.
\end{definition}

\begin{definition}[Eigenvector (EV) Method]
Let $A$ be an $n\times n$ PCM. The weight vector $w$ of $A$ determined by the EV method is defined as follows:
\begin{equation}
\label{eq:4}
    A\cdot w=\lambda_{\max}\cdot w,
\end{equation}
where the componentwise positive principal eigenvector $w$ is unique up to a scalar multiplication.
\end{definition}

These two methods are shown to be indeed similar in their results, however LLSM has significantly lower computational time \citep{Dong2008}.


In several situations some comparisons are absent, which may happen because the decision makers do not have time, willingness or possibility to make all of them, data have \color{black} been lost, the direct comparison is simply impossible (for instance in sports \citep{BozokiCsatoTemesi}), etc. When a PCM has missing elements, it is said to be an incomplete PCM (IPCM).

The LLSM and EV weight calculation methods can be generalized to the incomplete case as well, when the LLSM's optimization problem (Equation \ref{eq:3}) includes only the known elements of the matrix, while the EV method is based on the CR-minimal completion (CREV) of the PCM and its principal right eigenvector  \citep{Shiraishi1998,Shiraishi2002}. 

In this paper we analyze different kinds of filling in structures of IPCMs, thus we assume that the set of pairwise comparisons to be made can be chosen. We also heavily rely on the graph representation of IPCMs \citep{Gass1998}.

\begin{definition}[Representing graph of an IPCM]
An IPCM $A$ is represented by the undirected graph $G=(V,E)$, where the $V$ vertex set of $G$ corresponds to the alternatives (criteria) of $A$, and there is an edge in the edge set $E$ of $G$ if and only if the appropriate element of $A$ is known.
\end{definition}

We assume that no prior information is available about the items to be compared, thus in the examined filling in patterns we do not distinguish between the isomorphic representing graphs.
The optimal solutions of both above-mentioned weight calculation techniques for IPCMs (LLSM and CREV) are unique if and only if the representing graph is connected \citep{Bozoki2010}.

\begin{definition}[Connected graph]
In an undirected graph, two vertices $u$ and $v$ are called connected if the graph contains a path from $u$ to $v$. A graph is said to be connected if every pair of vertices in the graph is connected.
\end{definition}

The smallest connected systems are associated with spanning trees, which contain $n-1$ edges for $n$ vertices.

\begin{definition}[Spanning tree]
Let $G=(V,E)$ be a connected graph. $G^\prime=(V,E^\prime)$ is a spanning tree of $G$ if $E^\prime\subseteq E$ is a minimal set of edges that connect all vertices of $G$.
\end{definition}

An IPCM represented by a spanning tree can always be complemented to a consistent PCM, however, the results based on such an IPCM are usually extremely unreliable. The special importance of spanning trees is emphasized by the combinatorial weight calculation method \citep{Tsyganok2010}, which is built on the weight vectors obtained from all different spanning trees. This technique provides the same prioritization vector as the LLSM, if we use the geometric mean, both for PCMs \citep{Lundy2017} and IPCMs \citep{BozokiTsyganok2019}.

The results obtained by any weight calculation methods for IPCMs is strongly dependent on the number of known comparisons, namely the number of edges of the representing graph ($e$), and the arrangements of these known elements. Several properties have been examined in connection with the positioning of the known items, among which (some sense of) regularity of comparisons seems to be an especially important one \citep{WangTakahashi1998,Kulakowski,Szadoczki2020}, which can also be described by the representing graph.

\begin{definition}[$\boldsymbol{k}$-regularity]

A graph is called $k$-regular if every vertex has $k$ neighbours, which means that the degree of every vertex is $k$.
\end{definition}

When both the number of vertices ($n$) and the level of regularity ($k$) are odd, $k$-regularity is not possible. However, the graphs that are the closest to $k$-regularity in this case are called $k$-quasi-regular graphs \citep{Bozoki2020}.

\begin{definition}[$\boldsymbol{k}$-quasi-regularity]

A graph is called $k$-quasi-regular if exactly one vertex has degree $k+1$, and all the other vertices have degree $k$.

\end{definition}

In decision making the (quasi-)regularity of the representing graph ensures a certain level of symmetry, as every item is compared to the (approximately) same number of elements. This kind of property is also required in other fields, for instance, in the design of some sport tournaments \citep{Csato2017}.

We have only focused on multiplicative PCMs in the above definitions in this section, however, one can make the appropriate transformations to get an additive or a reciprocal (fuzzy) PCM from those \citep{brunelli2014introduction}. Thus, we would like to emphasize that all of our findings in the sections below are true for those types of matrices as well. 


\section{Methodology}
\label{sec:3}

Our aim is to find the filling structures that provide the closest results to the complete case for a given $(n,e)$ pair, number of alternatives (criteria) and comparisons. As it is assumed that we do not have any prior information, and so, the different items are not distinguished, we used Wolfram Mathematica \citep{Mathematica}, nauty and Traces \citep{McKay2014}, and IGraph/M \citep{szabolcshorvat2020} to generate every non-isomorphic (representing) graph for the examined $(n,e)$ pairs. Our extensive numerical simulations are based on the filling patterns related to these graphs.

The used methods are similar to \cite{Szadoczki2022}, however their study only focuses on a few special cases of filling patterns, their results are based on significantly smaller samples, and they even compare representing graphs with different number of edges (comparisons), i.e., cases where it is difficult to distinguish between the effect of the filling structure and the effect of the number of comparisons. While the current paper compares all the possible filling structures for a given $(n,e)$ pair (i.e., for the same number of comparisons) with a more general approach.

\color{black}

In order to measure the differences between the weight vectors, we apply commonly used cardinal and ordinal indicators, the Euclidean distance ($d_{euc}$) and the Kendall rank correlation coefficient (Kendall's $\tau$), respectively, which are defined as follows.

\begin{equation}
\label{eq:5}
d_{euc}(u,v)=\sqrt{\sum_{i=1}^{n}(u_i-v_i)^2}
\end{equation}

\begin{equation}
\label{eq:6}
\tau(u,v)=\frac{n_c(u,v)-n_d(u,v)}{n(n-1)/2}
\end{equation}

\noindent where $u$ denotes the weight vector obtained from a certain filling structure and $v$ is the weight vector computed from the complete PCM. $u$ and $v$ are normalized by $  \sum_{i=1}^{n} u_i = 1$, and $  \sum_{i=1}^{n} v_i = 1,$ respectively, and $v_i$ and $u_i$ denote the $i$th element of the appropriate vectors. $n_c(u,v)$ and $n_d(u,v)$ are the number of concordant and discordant pairs of the examined vectors, respectively. The range of the Kendall's $\tau$ is $[-1,1]$, and considering the notation in Equation \ref{eq:6}, a higher value indicates a better performance of the given filling pattern. However, in this case (Equation \ref{eq:5}) the Euclidean distance can be interpreted as an error, thus its smaller level  is preferred.

It is also worth mentioning that besides these, \cite{Szadoczki2022} used many different kinds of measures for the special cases examined by them, such as the Garuti index \citep{Garuti2020}, and all of those provided similar results, thus we chose to apply only the most well-known and basic metrics.

\color{black}

An instrumental part of our methodology is to determine the sample size needed in the simulations, which is based on a certain form of Chebyshev's inequality \citep{Steliga2010,Sawetal1984} that leads to the weak law of large numbers.

\begin{proposition}[Weak law of large numbers]
Let $(\xi_k)$ be independent and identically distributed random variables with finite standard deviation ($\sigma$), and let $E(.)$ denote the expected value operator. Then Equation \ref{eq:7} follows for all $\varepsilon > 0$:
\begin{equation}
\label{eq:7}
    P\left(\left|\frac{\sum_{k=1}^n \xi_k}{n}-E\left(\xi_1\right)\right| > \varepsilon\right) \leq \frac{\sigma^2}{n\varepsilon^2}\xrightarrow[n\to\infty]{} 0,
\end{equation}
where the last part of the expression means that the limit of the probability is $0$ as $n$ goes to infinity. 
\end{proposition}

The $\alpha=\sigma^2/n\varepsilon^2$ notation defines the significance level of our results, while $\varepsilon$ is the margin of error. We estimated the standard deviations of the Euclidean distances and the Kendall's $\tau$ measures for the different filling structures in our simulation and used an upper bound on it. Based on this method we applied a sample size of one million elements for every (representing) graph, which results in (as an upper bound as well)
\begin{itemize}
    \item $\alpha=0.01$ and $\varepsilon=0.0005$ for the computed Euclidean distances,
    \item and $\alpha=0.05$ and $\varepsilon=0.001$ for the calculated Kendall's $\tau$ measures.
\end{itemize}

As we mentioned earlier, the result of the EV weight calculation technique is similar to the LLSM, but its computational time is larger. This pattern is even stronger in the case of incompleteness (for CREV and incomplete LLSM, see for instance \cite{Csato2013}), thus due to the large sample sizes, in our simulations we mainly focus on the LLSM weight calculation technique. The results of the CREV method were computed for smaller cases ($n\leq5$) with a sample size of 500 000 as well, however, the ranking of filling patterns were always the same, and the indicators were almost always closer to the LLSM outcomes than the margin of error, thus we decided not to present them in much detail.

The process of the simulation for a given $(n,e)$ pair consisted of the following steps:

\begin{enumerate}
    \item $n$ random weights (in general they are denoted by $w_i$) were generated, where $w_i \in [1,9]$ is a uniformly distributed random real number for all $ i \in 1, 2, \ldots , n$. We calculated random $ n \times n$ complete and consistent PCMs, where the elements of the matrices were given by Equation \ref{eq:8}.
    \begin{equation}
    \label{eq:8}
a_{ij}=w_i/w_j
\end{equation}

\item Then three different perturbations of the items of consistent PCMs were used to get inconsistent matrices with three well-distinguishable inconsistency levels. These levels are denoted by weak, modest and strong given by Equations \ref{eq:9}, \ref{eq:10} and \ref{eq:11}.

\begin{equation}
\label{eq:9}
\hat{a}^{weak}_{ij}=\begin{cases} a_{ij}+\Delta_{ij} & :a_{ij}+\Delta_{ij}\geq 1\\ \frac{1}{1-\Delta_{ij}-(a_{ij}-1)}& :a_{ij}+\Delta_{ij}< 1 \end{cases} \hspace{2cm} \Delta_{ij} \in [-1,1]
\end{equation}

\begin{equation}
\label{eq:10}
\hat{a}^{modest}_{ij}=\begin{cases} a_{ij}+\Delta_{ij} & :a_{ij}+\Delta_{ij}\geq 1\\ \frac{1}{1-\Delta_{ij}-(a_{ij}-1)}& :a_{ij}+\Delta_{ij}< 1 \end{cases} \hspace{1.85cm} \Delta_{ij} \in \left[-\frac{3}{2},\frac{3}{2}\right]
\end{equation}

\begin{equation}
\label{eq:11}
\hat{a}^{strong}_{ij}=\begin{cases} a_{ij}+\Delta_{ij} & :a_{ij}+\Delta_{ij}\geq 1\\ \frac{1}{1-\Delta_{ij}-(a_{ij}-1)}& :a_{ij}+\Delta_{ij}< 1 \end{cases} \hspace{2cm} \Delta_{ij} \in [-2,2]
\end{equation}

Where $\hat{a}^{weak}_{ij}$, $\hat{a}^{modest}_{ij}$ and $\hat{a}^{strong}_{ij}$ are the elements of the perturbed PCMs, $a_{ij}$ is the element of the consistent PCM, $a_{ij}\geq1$ (we only perturb the elements above one and keep the reciprocity of the matrices), and $\Delta_{ij}$ is uniformly distributed in the given ranges. This perturbation method is able to produce ordinal differences as well (when $\hat{a}_{ij}<1$). It is  important to mention that  we account for the contrast that can be examined above and below $1$, thus our perturbed data is uniformly distributed around the original element on the scale presented by Figure~\ref{fig:1}, which also contains two examples. Our perturbation method aims to provide three different and meaningful inconsistency levels and it is, indeed, correlated with the Consistency Ratio (CR), as it is shown in Figure~\ref{fig:2}. We tested several combinations of parameters, and found that these resulted in the most relevant levels of CR, as the median CR value (based on a sample of 1000 matrices) is approximately 0.03, 0.065 and 0.1 for the weak, modest, and strong perturbations, respectively. We also ran the simulations with a modified perturbation method using lognormal errors in a multiplicative manner, however, they provided the exact same results and conclusions. For further details see Appendix~\ref{append:A}. \color{black} \\

\begin{center}
\begin{figure}[ht!]
\centering
  \includegraphics[width=0.8\textwidth]{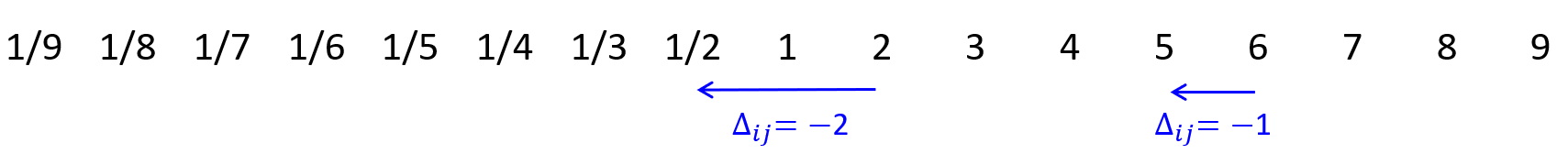}
\caption{The ratio scale $1/9, \ldots, 9$ and the perturbation of elements   according to (\ref{eq:9})--(\ref{eq:11}).}
\label{fig:1} 
\end{figure}
\end{center}

\input{CR2}

\bigskip

\item We deleted the respective elements of the matrices in order to get the filling structure that we were examining, and applied the LLSM (and CREV in  the  case of $n\leq5$) technique(s) to obtain the weights. The certain models' Euclidean distances and Kendall's $\tau$ measures were computed with respect to the weights that were calculated from the complete inconsistent matrices. The analyzed filling in patterns included all of those that can be represented by connected non-isomorphic graphs with parameters $(n,e)$.

\item We repeated steps 1-3 for $1 \ 000 \ 000$ times for every level of inconsistency (thus altogether we examined $3 \ 000 \ 000$ PCMs for a given $(n,e)$ pair). Finally, we saved the mean of Euclidean distances and Kendall's $\tau$ measures for the different filling in patterns.
\end{enumerate}

\begin{remark}
The distribution of the elements of complete PCMs is independent of $n$. This property holds for both consistent and perturbed complete PCM cases.
\end{remark}

  The reason behind this is that,  in the simulations at first the elements of a given matrix are generated independently from $n$, and then they are placed into the $n\times n$ PCM. The histograms of the complete PCM elements above $1$ in the different perturbation cases, based on samples containing 1 million elements each, are presented in Figure \ref{fig:3} (with a 0.1 bin width).

 \input{Histogram2}

According to the histograms, a higher level of perturbation (inconsistency) leads to a higher chance to have large (extreme) matrix elements.

\section{Results}
\label{sec:4}

We would like to emphasize that all of the results (and graph recommendations) presented in this section are under the following crucial assumptions.
\begin{enumerate}
    \item We can choose the comparisons that should be made (they are not given a priori).
    \item An ‘optimal' graph is the one that provides the closest LLSM weight vector on average to the one calculated from the complete matrix according to the measures presented in Section \ref{sec:3}.
    \item \label{enum:3} There is not any prior information about the items that should be compared, thus we can handle them in a symmetric way. This also means that the ‘reliability' and the weight \color{black} of every comparison is assumed to be the same.
\end{enumerate}
Naturally, if one or more of the assumptions above do not hold---for instance comparing to another benchmark instead of the complete inconsistent PCM---that could lead to other outcomes.

It is important to note that the interesting cases for our research  start  above three alternatives ($n$), as in  the  case of $n=2$ and $n=3$ there is always one non-isomorphic (representing) graph for every relevant pair of $(n,e)$ as it is shown in Figure \ref{fig:4}.

\input{n2_3}

The $n=4$ case also contains only a few possibilities, but it can be interesting in a decision problem, when there are several criteria and four alternatives, and it helps to understand the results for larger examples as well. Figure \ref{fig:5} presents the connected representing graphs for $n=4$ as a GRAPH of graphs.

\input{n4graph}

The value of $e$ is shown in every row of the GRAPH, in which an EDGE between two NODEs (=graphs) denotes that we can obtain one graph from the other one by adding (or deleting) exactly one edge.  
The GRAPH of graphs in Figure \ref{fig:5} is a 4-partite GRAPH with a further specific property, namely, that EDGEs go between levels $k$ and $k+1$ only ($k=1,2,3$). Note that if all EDGEs would be oriented ‘downwards' (i.e., the addition of an edge in the graph of comparisons), a partially ordered set of graphs (of comparisons) would be resulted in.
We denote the graph that provided the weight vectors with the smallest average Euclidean distance and the largest average Kendall's $\tau$ respect to the vectors calculated from the complete case by green background color for every $e$. If two optimal graphs are connected with an EDGE, then it is a partial optimal sequence, and the respective EDGE is also highlighted \color{black} by green. It is important to note that the relevant values for $e$ (the number of comparisons) are between $n-1$ (spanning trees) and $n(n-1)/2$ (complete graphs representing complete PCMs).

Among the spanning trees the star graph provided the smallest errors (Euclidean distances) and the largest Kendall's $\tau$ measures. This is not connected to the optimal graph with four edges, which is the $2$-regular cycle. However, from this point on the optimal graphs result in an optimal filling sequence. This is not surprising, as for $e=5$ and $e=6$ there is only one possible non-isomorphic representing graph, but this example probably helps to understand the following cases. Tables \ref{tab:1} and \ref{tab:2} present the results provided by the graphs with $(n=4,e=3)$ and $(n=4,e=4)$ respectively, in  the  case of the different perturbation levels. The name of the optimal graph, and the best values in every column are highlighted by bold text. \color{black}

\begin{table}[ht!] \centering

\begin{tabular}{@{}ccccccc@{}}\toprule
	Graph & \multicolumn{2}{c}{Weak} & \multicolumn{2}{c}{Modest} & \multicolumn{2}{c}{Strong} \\ \hline
	 & $d_{euc}$ & Kendall's $\tau$ & $d_{euc}$ & Kendall's $\tau$ & $d_{euc}$ & Kendall's $\tau$ \\ \hline
	\textbf{Star graph} & \textbf{0.0918} & \textbf{0.7306} & \textbf{0.1293} & \textbf{0.6639} & \textbf{0.1620} & \textbf{0.6164} \\
	Line graph &0.0967 & 0.7194 & 0.1361 & 0.6501 & 0.1701 & 0.6020  \\
	\bottomrule
	\end{tabular}
 \caption{The average Euclidean distances and Kendall's $\tau$ measures for the graphs with $(n=4,e=3)$ in  the  case of the different perturbation levels.}
 \label{tab:1}
\end{table}

\begin{table}[H] \centering

\begin{tabular}{@{}ccccccc@{}}\toprule
	 Graph & \multicolumn{2}{c}{Weak} & \multicolumn{2}{c}{Modest} & \multicolumn{2}{c}{Strong} \\ \hline
	 & $d_{euc}$ & Kendall's $\tau$ & $d_{euc}$ & Kendall's $\tau$ & $d_{euc}$ & Kendall's $\tau$ \\ \hline
	Not regular graph & 0.0650 & 0.8027 & 0.0920 & 0.7496 & 0.1156 & 0.7111 \\
	$\boldsymbol{2}$\textbf{-regular graph} & \textbf{0.0543} &  \textbf{0.8216} &  \textbf{0.0771} &  \textbf{0.7705} & \textbf{0.0970} & \textbf{0.7328}  \\
	\bottomrule
	\end{tabular}
 \caption{The average Euclidean distances and Kendall's $\tau$ measures for the graphs with $(n=4,e=4)$ in  the case of the different perturbation levels.}
 \label{tab:2}
\end{table}

Based on Tables \ref{tab:1} and \ref{tab:2}, one can observe that for a given $(n,e)$ pair,  the same graphs provided the best results on average for the examined measures for every perturbation level. There are indeed significant differences between the examined graphs (the margin of error is $0.0005$ for the Euclidean distances and $0.001$ for the Kendall's $\tau$ measures). It is also easy to see that a stronger perturbation results in higher errors, while an additional edge leads to smaller distances and higher ordinal correlations. Figure \ref{fig:6} presents the relation between the number of comparisons ($e$) and the analyzed cardinal ($d_{euc}$) and ordinal (Kendall's $\tau$) measures, which can help practitioners to determine the minimal sufficient number of comparisons in a given problem. Note that Figure \ref{fig:6} shows the results for the optimal graphs for every $e$, thus one optimal value is not necessarily reachable from the previous one, only in  the  case of partial optimal sequences.

 If we know in advance that the decision maker is willing to provide exactly $e = n-1 = 3$ comparisons, then, according to Figure \ref{fig:5}, we recommend the star graph, i.e., filling in one (e.g. the first) row/column of the pairwise comparison matrix, namely elements $a_{12}, a_{13}$ and $a_{14}$ (in any order), also summarized in Table \ref{tab:3}.  It is worth noting that in this case Assumption \ref{enum:3} has a special importance as all items are compared to one (pivotal) item.

\begin{table}[htbp!]
\centering
\begin{tabular}{|c|c|c|c|c|}
\hline
   &  1  &  2  &  3  &  4    \\ \hline
1  &     & \#1' & \#2' & \#3'    \\ \hline
2  &     &     &     &       \\ \hline
3  &     &     &     &       \\ \hline
4  &     &     &     &       \\     \hline
\end{tabular}
\caption{Filling in sequence for $n = 4$, $e = 3$. Orders with ' are interchangeable.} 
\label{tab:3}
\end{table}

 If we assume that the decision maker is willing to provide more than three comparisons, the optimal filling in sequence is $
\{a_{12}, 
a_{23}, 
a_{34}, 
a_{14}\}$  (the first four  comparisons can be made in any order), 
followed by
$a_{13}$ and finally 
$a_{24},$
also summarized in Table \ref{tab:4}.

\begin{table}[H]
\centering
\begin{tabular}{|c|c|c|c|c|}
\hline
   &  1  &  2   &  3   &  4      \\ \hline
1  &     & \#1' & \#5  & \#4'         \\ \hline
2  &     &      & \#2' & \#6         \\ \hline
3  &     &      &      & \#3'    \\ \hline
4  &     &      &      &         \\     \hline
\end{tabular}
\caption{Filling in sequence for $n = 4$, $e > 3$. Orders with ' are interchangeable.} 
\label{tab:4}
\end{table}

\input{n4Optimal}

Naturally, in a general practical application, we do not know how many comparisons will be provided by the decision maker. However, if the number of comparisons is lower than $e=3(=n-1)$, then it is not possible to determine a unique weight vector, while if it is strictly larger than $e=3$, then fortunately it happens that a partial optimal sequence leads to the complete matrix. Based on that we would consider the instance with exactly $e=3$ comparisons a special case, and follow the partial optimal filling in sequence otherwise. \color{black}

For larger number of alternatives (criteria, $n$), the possible number of connected graphs increases quickly, thus it is even more relevant to determine the optimal filling structure. In  the  case of $n=5$, there are $21$ connected graphs altogether. Their 7-partite GRAPH of graphs can be seen in Figure \ref{fig:7}, using the same notations as before.

\input{n5graph}

One can see many similarities with the previous outcomes. The star graph resulted in the smallest Euclidean distance and the largest Kendall's $\tau$ measure among the spanning trees, once again. It is not connected to the optimal graph with $e=5$, which is the $2$-regular cycle, as before. The next optimal graph with $e=6$ is not connected to the cycle, as well, however, from that point on there is a partial optimal sequence to the complete filling of the represented PCM. Somewhat surprisingly, the graphs providing the smallest Euclidean distances resulted in the largest Kendall's $\tau$ for every single case, except for $e=8$. However, in that case the difference between the Kendall's $\tau$ measures for the two possible graphs is within the margin of error, thus we highlighted the graph that is better according to the Euclidean distance, and it \color{black} is better in both indicators if we consider the CREV weight calculation technique or if we use lognormal perturbations (see Appendix \ref{append:A})\color{black}. It is worth mentioning that this graph is the $3$-quasi-regular graph on $n=5$.

Figure \ref{fig:8} shows the relation between the number of comparisons ($e$) and the analyzed measures for $n=5$ in  the  case of optimal graphs. One optimal value is not necessarily reachable from the previous one, as before. Minimal thresholds could be determined for the number of comparisons based on this figure for certain decision problems.

If we know in advance that the decision maker is willing to provide exactly $e = n-1 = 4$ comparisons, then, according to Figure \ref{fig:7}, we recommend the star graph  (with special attention to Assumption \ref{enum:3}), i.e., filling in one (e.g. the first) row/column of the pairwise comparison matrix, namely elements $a_{12}, a_{13}, a_{14}$ and $a_{15}$, also summarized in Table \ref{tab:5}.

\begin{table}[ht!]
\centering
\begin{tabular}{|c|c|c|c|c|c|}
\hline
   &  1  &  2   &  3   &    4   &    5    \\ \hline
1  &     & \#1' & \#2' &   \#3' &   \#4'  \\ \hline
2  &     &      &      &        &         \\ \hline
3  &     &      &      &        &         \\ \hline
4  &     &      &      &        &         \\ \hline
5  &     &      &      &        &         \\     \hline
\end{tabular}
\caption{Filling in sequence for $n = 5$, $e = 4$. Orders with ' are interchangeable.} 
\label{tab:5}
\end{table}

In  the  case the decision maker is willing to provide exactly $e = n = 5$ comparisons, then we should  make the comparisons  along an $n$-cycle, e.g.,
$\{a_{12}, 
a_{23}, 
a_{34}, 
a_{45}$ and 
$a_{15}\}$ (the five  comparisons can be made  in any order), 
 also summarized in Table \ref{tab:6}.

\begin{table}[ht!]
\centering
\begin{tabular}{|c|c|c|c|c|c|}
\hline
   &  1  &  2   &  3   &    4   &    5    \\ \hline
1  &     & \#1' &      &        &  \#5'   \\ \hline
2  &     &      & \#2' &        &         \\ \hline
3  &     &      &      &  \#3'  &         \\ \hline
4  &     &      &      &        &  \#4'   \\ \hline
5  &     &      &      &        &         \\     \hline
\end{tabular}
\caption{Filling in sequence for $n = 5$, $e = 5$. Orders with ' are interchangeable.} 
\label{tab:6}
\end{table}

When we can assume that the decision maker is willing to provide more than five comparisons, the optimal filling in sequence is $
\{
a_{14},
a_{15},
a_{24},
a_{25}, 
a_{34}, 
a_{35}\}$ (the first six  comparisons can be made  in any order), 
followed by $
a_{12},
a_{13},
a_{45},
$
and finally 
$a_{23},$  also summarized in Table \ref{tab:7}.

\begin{table}[ht!]
\centering
\begin{tabular}{|c|c|c|c|c|c|}
\hline
   &  1  &   2   &  3    &    4   &    5    \\ \hline
1  &     &  \#7     & \#8      &  \#1'  &  \#2'   \\ \hline
2  &     &       &   \#10    &  \#3'  &  \#4'   \\ \hline
3  &     &       &       &  \#5'  &  \#6'   \\ \hline
4  &     &       &       &        &    \#9     \\ \hline
5  &     &       &       &        &         \\ \hline
\end{tabular}
\caption{Filling in sequence for $n = 5$, $e > 5$. Orders with ' are interchangeable.} 
\label{tab:7}
\end{table}

Again, in an application, it is difficult to know the number of provided comparisons in advance. Still, if the decision maker provides at least $e=6(=n+1)$ comparisons, then there is a partial optimal filling in sequence to the complete PCM, while if the number of provided comparisons is at most $e=3(=n-2)$, then it is not possible to determine a unique weight vector. Thus, the instances with exactly $e=4(=n-1)$ and $e=5(=n)$ comparisons are considered to be special cases, and generally the partial optimal filling in sequence is recommended to be followed. \color{black}

\input{n5Optimal}

Finally, for $n=6$, there are $112$ possible connected (representing) graphs. Figure \ref{fig:9} shows the 11-partite GRAPH of graphs for this case, however, in order to keep it visible, we only denote each possible graph by a vertex, and present the optimal cases in detail in Figure \ref{fig:10}. For $e=12$ and $e=13$ the results are close to each other, and some of the differences of the Kendall's $\tau$ measures are also smaller than the margin of error. Here the best graph according to the Euclidean distance and the Kendall's $\tau$ are different as well. However, we highlighted the graphs which were at least second according to at least one indicator by a lighter green color. These highlighted graphs for a given $e$ practically provide the same results. As there is always a unique optimal graph according to the Euclidean distance, we denoted those with an $E$. We have not highlighted the EDGEs by green color on this part of the GRAPH of graphs, because of the similar results (ties), however, the EDGE connecting the graphs (NODEs) that provided the best results according to the Euclidean distance are highlighted by blue color. In Figure \ref{fig:10} for $e=12$ and $13$ also the graphs, which are the best according to the Euclidean distance are presented.

 If we know in advance that the decision maker is willing to provide exactly $e = n-1 = 5$ comparisons, then, according to Figure \ref{fig:9}, we recommend the star graph  (with special attention to Assumption \ref{enum:3}), i.e., filling in one (e.g. the first) row/column of the pairwise comparison matrix, namely elements $a_{12}, a_{13}, a_{14}, a_{15}$ and $a_{16}$, 
 also summarized in Table \ref{tab:8}.

\begin{table}[ht!]
\centering
\begin{tabular}{|c|c|c|c|c|c|c|}
\hline
   &  1  &  2   &  3   &    4   &    5    &    6    \\ \hline
1  &     & \#1' & \#2' &   \#3' &   \#4'  &  \#5'   \\ \hline
2  &     &      &      &        &         &         \\ \hline
3  &     &      &      &        &         &         \\ \hline
4  &     &      &      &        &         &         \\ \hline
5  &     &      &      &        &         &         \\ \hline
6  &     &      &      &        &         &         \\     \hline
\end{tabular}
\caption{Filling in sequence for $n = 6$, $e = 5$. Orders with ' are interchangeable.} 
\label{tab:8}
\end{table}

 If the decision maker  is willing to provide more than five comparisons, the recommended filling in sequence is  $
\{
a_{14},
a_{15},
a_{24},
a_{26}, 
a_{35}, 
a_{36}\}$ (the first six  comparisons can be made  in any order), 
followed by $
\color{orange} a_{25},    
a_{34}, \color{black}       
a_{16},                   
a_{12},                   
a_{46},                   
\color{orange} a_{23},    
a_{45},  \color{black}      
a_{56},                   
$
and finally 
$a_{13},$                 
also summarized in Table \ref{tab:9}.

\begin{table}[ht!]
\centering
\begin{tabular}{|c|c|c|c|c|c|c|}
\hline
   &  1  &  2   &  3   &    4   &    5    &    6    \\ \hline
1  &     &  \#10    &   \#15   &  \#1'  &  \#2'   &   \#9   \\ \hline
2  &     &      &   \color{orange}\#12   &  \#3'  &   \color{orange}\#7   &  \#4'   \\ \hline
3  &     &      &      &  \color{orange}\#8   &  \#5'   &  \#6'   \\ \hline
4  &     &      &      &        &   \color{orange} \#13     &  \#11   \\ \hline
5  &     &      &      &        &         &   \#14      \\ \hline
6  &     &      &      &        &         &         \\     \hline
\end{tabular}
\caption{Filling in sequence for $n = 6$, $e > 5$. Orders with ' are interchangeable.} 
\label{tab:9}
\end{table}

Since there is no path along all the optimal graphs, the filling in sequence above includes as many as possible. The remaining EDGEs are colored with orange in Figure \ref{fig:9}, and it should be noted \color{black} that the other included graphs are as close to optimal ones as possible. In general it is difficult to predict the number of comparisons provided by the decision maker, thus we recommend to follow this sequence of comparisons, and treat the instance, when $e=5(=n-1)$ as a special case. \color{black}

\input{n6graph}

\input{MKE_PTE}

One can observe many similarities with the earlier outcomes in connection with the concrete graphs, and the pattern of optimal graphs as well. Among the spanning trees, the star graph provided the best results according to both measures again. For $e=6$ the $2$-regular cycle turned out to be the optimal case, just as earlier. The optimal graphs with $e=5,6$ and $7$ are not connected, but from that point on we can determine an optimal filling in sequence to the complete graph (if we consider all the light green cases optimal).

Moreover, for $e=9$ the optimal graph is the single bipartite $3$-regular graph on six vertices, while for $e=12$ the highlighted graph, which provided the best results according to the Euclidean distance and the second best according to the Kendall's $\tau$, is the only $4$-regular graph on six vertices. Based on the general conclusions of the \color{black} simulations, we can make several important remarks.

\begin{remark}
The star graph provided the best results according to both measures for all examined $(n,e=n-1)$ cases. Thus we can say that it is an optimal structure, intuitively it keeps this property for larger cases $(n)$, as well.
\end{remark}

\begin{remark}
For the $(n,e=k\cdot n/2)$ examples, the optimal graph is always a $k$-regular graph. Furthermore, $k$-quasi-regular graphs are optimal as well. One can say that regularity is indeed important in a more general way, as in all of the examined instances, the degree of different vertices (the number of comparisons) are as close as possible.
\end{remark}

\begin{remark}
The optimal graphs are always bipartite graphs, or the closest ones to that.
\end{remark}

The analyzed indicators for optimal graphs in the case of different number of comparisons ($e$) can be seen in Figure \ref{fig:11} for $n=6$. Again, it can serve as a guide for practitioners.

\input{n6Optimal}

All of our simulation results provided optimal filling structures (representing graphs) for the examined $(n,e)$ pairs, as well as (partial) optimal filling sequences. The outcomes show indeed similar patterns for different parameters, and can support both applications and theoretical studies.

Finally, it is worth mentioning that although a practical MCDM problem usually has several hierarchical levels,---thus many PCMs have to be filled in and the overall number of comparisons is high,---the size of the matrices usually do not exceed $6\times6$ \citep{Abele2018}, which makes the presented results even more relevant. As shown in Figure \ref{fig:4}, finding the optimal sequence of comparisons is only interesting for at least $4$ alternatives, while for $7$ or more alternatives the problem becomes computationally too demanding to proceed. There are $112$ possible connected (representing) graphs for $n=6$ ($112$ NODEs in Figure \ref{fig:9}), while the same number for $n=7,8,$ and $9$ are $853,$ $11117,$ and $261080$, respectively. At the same time, the difference between the different filling patterns seems to be smaller as $n$ increases. Thus, there is a good chance that there would be many cases, where the difference between the graphs would be smaller than the margin of error for larger $n$ values. Furthermore, based on our research, \cite{Gyarmati2023} found that the exact same graphs are optimal in the case of fundamentally different models as well, which are based on paired comparisons, i.e., the Bradley--Terry and Thurstone models. These findings suggest that our results are rather general and not model-specific.

\section{Conclusion and further research}
\label{sec:5}

In this paper we analyzed all possible filling structures of incomplete pairwise comparison matrices when there is no prior information available for the compared items, in the case of at most six alternatives (criteria). The study heavily relied on the representing graphs of pairwise comparisons as well as on extensive numerical simulations with large samples. We compared the weight vectors (calculated by the incomplete LLSM) related to the certain filling patterns and compared them based on their Euclidean distance and Kendall's $\tau$ measure with the weights obtained from the complete case.

We found that there is a strong connection between the examined cardinal and ordinal indicators,  thus we could determine the best filling structure for a given number of alternatives and comparisons, which is a significant finding itself. However, one of the main contributions of the paper is that, many of these optimal graphs resulted in optimal filling in sequences as illustrated by different paths in the examined GRAPHs of graphs first in the literature. \color{black}

The filling structure represented by a star graph turned out to be optimal among the graphs (filling patterns) with the same cardinality (spanning trees). Regular graphs also seem to provide optimal solutions, and regularity is a common property of the optimal cases in a more general sense. 

Both theorists and practitioners can utilize our findings not just to apply the optimal filling structure in their problems, but also to use the optimal filling sequences in decision making problems where the decision maker can abandon the problem at any period of the process. Furthermore, our results on the difference between the optimal patterns and the complete case for different number of comparisons can serve as a guide to determine the minimal sufficient number of comparisons for a given problem.

The presented results seem to be robust regarding
\begin{itemize}
    \item the weight calculation technique (eigenvector or logarithmic least squares method);
    \item the level of inconsistency (weak, modest and strong perturbation levels);
    \item the way of perturbation (multiplicative lognormal errors or uniformly distributed errors on a modified scale);
    \item the used distance metrics (Euclidean distance and Kendall's $\tau$);
    \item even the model using pairwise comparisons, considering the findigs of \cite{Gyarmati2023}.
\end{itemize}

A future research can investigate the certain comparisons that decrease the errors the most during the filling in process. When should we stop to ask even more questions from the decision maker? Do the last few comparisons provide significant information? How does this problem relate to the representing graph?

Empirical PCMs may differ from simulated ones, and many collections of those matrices (even with the complete filling in order) are available \citep{Bozoki2013}, thus in a future research it is important to test our findings on empirical matrices as well.

Naturally larger cases, other weight calculation methods and different distance measures can be further investigated as well. Are the findings remain true for a large number of alternatives? How much are they dependent on the used techniques and measures? What can we say when some prior information, for instance, the best or the worst alternatives, perhaps both, are known \citep{Rezaei2015,Mustajoki2005,Edwards1994,Winterfeldt1986}?

Our results can be useful in other areas as well, for instance, in designing sport tournaments. If we would like to plan the different rounds, we should make a number of comparisons simultaneously. This leads to the general question: besides optimal direct sequences, how does the optimal graphs include each other (indirectly)?

\section*{Acknowledgements}
The authors thank the valuable comments and suggestions of the anonymous Reviewers. The authors are grateful for the comments of László Csató and Éva Orbán-Mihálykó. The project identified by EKOP-CORVINUS-24-4-080 was realized with the support of the National Research, Development, and Innovation Fund provided by the Ministry of Culture and Innovation, as part of the University Research Scholarship Program announced for the 2024/2025 academic year. The research was supported by the National Research, Development and Innovation Office under Grants FK 145838 and TKP2021-NKTA-01 NRDIO.

\section*{Declaration}

\subsection*{Competing Interest}
The authors have no relevant financial or non-financial interests to disclose.


\bibliographystyle{apalike} 
\bibliography{main}
\addcontentsline{toc}{section}{References}

\newpage

\newcounter{appendfigure}
\setcounter{appendfigure}{1}

\newcounter{append}
\renewcommand{\theappend}{A}

\section*{Appendix A}
\addcontentsline{toc}{section}{Appendix A}
\refstepcounter{append}
\label{append:A}

A modified version of the perturbations---using lognormal errors in a multiplicative manner---was also investigated, where the weak, modest and strong inconsistency levels were ensured by Equations \eqref{appendequation:1}, \eqref{appendequation:2}, and \eqref{appendequation:3}.

\begin{equation}
\label{appendequation:1} \tag{A1}
\hat{a}^{weak}_{ij}= a_{ij}\cdot\Delta_{ij} \hspace{2cm} \Delta_{ij} \sim Lognorm(0,0.35)
\end{equation}

\begin{equation}
\label{appendequation:2} \tag{A2}
\hat{a}^{modest}_{ij}= a_{ij}\cdot\Delta_{ij} \hspace{2cm} \Delta_{ij} \sim Lognorm(0,0.5)
\end{equation}

\begin{equation}
\label{appendequation:3} \tag{A3}
\hat{a}^{strong}_{ij}= a_{ij}\cdot\Delta_{ij} \hspace{2cm} \Delta_{ij} \sim Lognorm(0,0.625)
\end{equation}

Where the same notation is used as in the main text, $\hat{a}^{weak}_{ij}$, $\hat{a}^{modest}_{ij}$ and $\hat{a}^{strong}_{ij}$ are the elements of the perturbed PCMs, $a_{ij}$ is the element of the consistent PCM, $a_{ij}\geq1$ (we only perturb the elements above one and keep the reciprocity of the matrices), and $\Delta_{ij}$ is a lognormal error with given parameters (expected value and standard deviation of the underlying normal distribution). These perturbations result in the inconsistency levels presented by Figure \ref{fig:append1} via Box plots.

\input{CR_lognorm}

This version of the perturbations basically provide the same levels of inconsistency as the one presented in Figure \ref{fig:2}, except that the maximum values are tend to be more extreme in this case. This is caused by the fact that there is no theoretical upper bound on the maximum of the CR in this case (no theoretical bounds on the perturbation) compared to the uniform case, however, the probability of these extreme cases are quite small.

We ran all of the simulations (using the LLSM weight calculation technique) with these perturbations as well, and they provided exactly the same results and conclusions as the ones presented in Section \ref{sec:4}.

\end{document}

%% file: CR2.tex
\begin{figure}[ht!]
	\centering
	\pgfplotsset{scaled y ticks=false}
	\pgfplotsset{every axis plot/.append style={line width=1.5pt}}
	\begin{tikzpicture}
	\begin{scope}
\begin{axis}[box plot width=1.5mm,
title = {Weak perturbation},
ylabel = Consistency Ratio (CR),
ymin=0,ymax=0.4,
xtick={5,6,7,8,9,10},
width=5cm,
height=8cm]
\boxplot[forget plot, blue,box plot whisker bottom index=1,
    box plot whisker top index=5,
    box plot box bottom index=2,
    box plot box top index=4,
    box plot median index=3]{weak.dat};

\end{axis}
\end{scope}

\begin{scope}[xshift=5.5cm]
\begin{axis}[box plot width=1.5mm,
ymin=0,ymax=0.4,
xtick={5,6,7,8,9,10},
title = {Modest perturbation},
xlabel = Number of alternatives/criteria ($n$),
width=5cm,
height=8cm]

\boxplot[forget plot, red,box plot whisker bottom index=1,
    box plot whisker top index=5,
    box plot box bottom index=2,
    box plot box top index=4,
    box plot median index=3]{modest.dat};

\end{axis}
\end{scope}

\begin{scope}[xshift=11cm]
\begin{axis}[box plot width=1.5mm,
ymin=0,ymax=0.4,
xtick={5,6,7,8,9,10},
title = {Strong perturbation},
width=5cm,
height=8cm]

\boxplot[forget plot, green,box plot whisker bottom index=1,
    box plot whisker top index=5,
    box plot box bottom index=2,
    box plot box top index=4,
    box plot median index=3]{strong.dat};

\end{axis}
\end{scope}

\end{tikzpicture}

\caption{The relation between CR and our element-wise perturbation via Box plots. Each Box plot is based on 1000 randomly generated perturbed PCMs, and shows the minimum, maximum, and the first, second (median), and third quartile of the gained CR values.}
\label{fig:2}
\end{figure}

%% file: Histogram2.tex
\begin{figure}[H]
\centering
	\pgfplotsset{scaled y ticks=false}
	\pgfplotsset{every axis plot/.append style={line width=1.5pt}}

\begin{tikzpicture}
\begin{axis}[
    width=15cm,
    xlabel={Value of elements},
    ylabel={Number of elements},
    xmin=1, xmax=9,
    ymin=0, ymax=70000,
    xtick={1,3,5,7,9},
    legend pos=north west,
    ymajorgrids=true,
    grid style=dashed,
    legend style ={ at={(0.63,0.95)},
        anchor=north west, 
        fill=white,align=left},
]

\addplot[
    color=orange,
    ]
    coordinates {
    (1,113914)	(1.1,94122)	(1.2,79528)	(1.3,68263)	(1.4,58531)	(1.5,51685)	(1.6,44905)	(1.7,39590)	(1.8,35422)	(1.9,31419)	(2,28652)	(2.1,25950)	(2.2,23585)	(2.3,21334)	(2.4,19367)	(2.5,18068)	(2.6,16436)	(2.7,15201)	(2.8,14038)	(2.9,12889)	(3,12106)	(3.1,11104)	(3.2,10361)	(3.3,9809)	(3.4,9077)	(3.5,8358)	(3.6,7983)	(3.7,7475)	(3.8,6947)	(3.9,6513)	(4,6094)	(4.1,5774)	(4.2,5395)	(4.3,4981)	(4.4,4742)	(4.5,4474)	(4.6,4177)	(4.7,4104)	(4.8,3859)	(4.9,3641)	(5,3428)	(5.1,3126)	(5.2,3058)	(5.3,2846)	(5.4,2675)	(5.5,2544)	(5.6,2481)	(5.7,2221)	(5.8,2137)	(5.9,1992)	(6,1914)	(6.1,1800)	(6.2,1627)	(6.3,1554)	(6.4,1483)	(6.5,1407)	(6.6,1284)	(6.7,1235)	(6.8,1142)	(6.9,1071)	(7,1005)	(7.1,974)	(7.2,853)	(7.3,813)	(7.4,709)	(7.5,715)	(7.6,616)	(7.7,543)	(7.8,495)	(7.9,451)	(8,392)	(8.1,328)	(8.2,295)	(8.3,252)	(8.4,201)	(8.5,180)	(8.6,114)	(8.7,102)	(8.8,47)	(8.9,17)

    };
    \addlegendentry{Consistent}
 
 \addplot[
    color=blue, dashed,
    ]
    coordinates {
    (1,61327)	(1.1,61504)	(1.2,60590)	(1.3,59439)	(1.4,58394)	(1.5,56428)	(1.6,54121)	(1.7,51554)	(1.8,47946)	(1.9,43426)	(2,38276)	(2.1,33254)	(2.2,29667)	(2.3,26340)	(2.4,24036)	(2.5,21217)	(2.6,19375)	(2.7,17648)	(2.8,16134)	(2.9,14870)	(3,13780)	(3.1,12513)	(3.2,11681)	(3.3,10855)	(3.4,9977)	(3.5,9314)	(3.6,8701)	(3.7,7938)	(3.8,7633)	(3.9,7167)	(4,6771)	(4.1,6266)	(4.2,5918)	(4.3,5601)	(4.4,5140)	(4.5,4911)	(4.6,4514)	(4.7,4261)	(4.8,4077)	(4.9,3739)	(5,3550)	(5.1,3404)	(5.2,3197)	(5.3,3035)	(5.4,2855)	(5.5,2695)	(5.6,2532)	(5.7,2463)	(5.8,2163)	(5.9,2066)	(6,2043)	(6.1,1837)	(6.2,1740)	(6.3,1685)	(6.4,1660)	(6.5,1443)	(6.6,1416)	(6.7,1246)	(6.8,1187)	(6.9,1160)	(7,1037)	(7.1,1052)	(7.2,919)	(7.3,859)	(7.4,717)	(7.5,733)	(7.6,637)	(7.7,558)	(7.8,536)	(7.9,453)	(8,444)	(8.1,386)	(8.2,349)	(8.3,273)	(8.4,237)	(8.5,233)	(8.6,177)	(8.7,148)	(8.8,134)	(8.9,110)

    };
    \addlegendentry{Weak perturbation}

    \addplot[
    color=red, dotted,
    ]
    coordinates {
(1,48811)	(1.1,48906)	(1.2,48412)	(1.3,48509)	(1.4,47896)	(1.5,47311)	(1.6,46330)	(1.7,45812)	(1.8,44728)	(1.9,43239)	(2,41603)	(2.1,40022)	(2.2,37964)	(2.3,35107)	(2.4,31960)	(2.5,28411)	(2.6,25023)	(2.7,22369)	(2.8,19936)	(2.9,18046)	(3,16228)	(3.1,15046)	(3.2,13645)	(3.3,12555)	(3.4,11501)	(3.5,10885)	(3.6,9872)	(3.7,9034)	(3.8,8450)	(3.9,7814)	(4,7461)	(4.1,6843)	(4.2,6464)	(4.3,6123)	(4.4,5706)	(4.5,5250)	(4.6,5130)	(4.7,4744)	(4.8,4327)	(4.9,4186)	(5,3837)	(5.1,3688)	(5.2,3393)	(5.3,3371)	(5.4,3165)	(5.5,2823)	(5.6,2828)	(5.7,2582)	(5.8,2464)	(5.9,2236)	(6,2124)	(6.1,2056)	(6.2,1881)	(6.3,1754)	(6.4,1651)	(6.5,1602)	(6.6,1466)	(6.7,1339)	(6.8,1238)	(6.9,1145)	(7,1069)	(7.1,988)	(7.2,966)	(7.3,885)	(7.4,773)	(7.5,779)	(7.6,700)	(7.7,614)	(7.8,579)	(7.9,498)	(8,456)	(8.1,412)	(8.2,396)	(8.3,362)	(8.4,329)	(8.5,296)	(8.6,259)	(8.7,208)	(8.8,200)	(8.9,186)

    };
    \addlegendentry{Modest perturbation}

    \addplot[
    color=green, mark = star,
    ]
    coordinates {
    (1,40987)	(1.1,40741)	(1.2,40326)	(1.3,39928)	(1.4,40140)	(1.5,39832)	(1.6,39597)	(1.7,38965)	(1.8,38482)	(1.9,38246)	(2,38021)	(2.1,37242)	(2.2,36517)	(2.3,35491)	(2.4,33968)	(2.5,33139)	(2.6,31639)	(2.7,29525)	(2.8,27529)	(2.9,25064)	(3,22250)	(3.1,19805)	(3.2,17683)	(3.3,16061)	(3.4,14417)	(3.5,13065)	(3.6,11780)	(3.7,10972)	(3.8,10237)	(3.9,9397)	(4,8694)	(4.1,8069)	(4.2,7532)	(4.3,7030)	(4.4,6398)	(4.5,5940)	(4.6,5571)	(4.7,5256)	(4.8,4841)	(4.9,4639)	(5,4330)	(5.1,4056)	(5.2,3699)	(5.3,3564)	(5.4,3375)	(5.5,3170)	(5.6,2881)	(5.7,2776)	(5.8,2626)	(5.9,2445)	(6,2239)	(6.1,2186)	(6.2,2108)	(6.3,2006)	(6.4,1822)	(6.5,1747)	(6.6,1564)	(6.7,1390)	(6.8,1433)	(6.9,1326)	(7,1176)	(7.1,1099)	(7.2,1063)	(7.3,967)	(7.4,935)	(7.5,898)	(7.6,768)	(7.7,731)	(7.8,617)	(7.9,640)	(8,580)	(8.1,551)	(8.2,477)	(8.3,418)	(8.4,405)	(8.5,362)	(8.6,352)	(8.7,333)	(8.8,265)	(8.9,231)

    };
    \addlegendentry{Strong perturbation}

\end{axis}

\end{tikzpicture}

\caption{The histograms of the $\geq 1$ elements of PCMs in case of different perturbations based on a sample of 1 million elements (with a 0.1 bin width).}
\label{fig:3}
\end{figure}

%% file: n2_3.tex
\begin{figure}[ht!]
\centering
\begin{tikzpicture}[every node/.style={circle,inner sep=2pt,draw=black,fill=black}]

\tikzstyle{block} = [circle,inner sep=2pt,draw=white,fill=white];

 \node (21) at (-9,-1.5) {};
  \node (22) at (-8,-1.5) {};
  
  \node [block] at (-8.5,0) (2) {$n=2, e=1$};
  
  \node (31) at (-4,-1.5) {};
  \node (32) at (-5,-2.5) {};
  \node (33) at (-3,-2.5) {};
  
  \node [block] at (-4,0) (2) {$n=3, e=2$};
  
  \node (311) at (1,-1.5) {};
  \node (322) at (2,-2.5) {};
  \node (333) at (0,-2.5) {};
  
  \node [block] at (1,0) (2) {$n=3, e=3$};
  
  \draw [draw=black!50] (21) -- (22)
                        (31) -- (32)
                        (31) -- (33)
                        (311) -- (322)
                        (311) -- (333)
                        (322) -- (333);

\end{tikzpicture}

\caption{The connected non-isomorphic representing graphs for $n\leq3$.}
\label{fig:4}
\end{figure}

%% file: n4graph.tex
\begin{figure}[H]
\begin{center}
\begin{tikzpicture}[every node/.style={circle,inner sep=2pt,draw=black,fill=black!20}]

\tikzstyle{node2} = [draw=white,fill=white];

 \node (11) at (0,0) {};
  \node (12) at (0,-1) {};
  \node (13) at (-1,-1) {};
  \node (14) at (1,-1) {};
  
  \node (21) at (5,0) {};
  \node (22) at (6,0) {};
  \node (23) at (5,-1) {};
  \node (24) at (6,-1) {};

  \node [node2] (3el) at (-5,0) {\text{ $e=3$ }};
  
  \node (31) at (-0.5,-4) {};
  \node (32) at (0.5,-4) {};
  \node (33) at (-0.5,-5) {};
  \node (34) at (0.5,-5) {};
  
  \node (41) at (5,-4) {};
  \node (42) at (6,-4) {};
  \node (43) at (5,-5) {};
  \node (44) at (6,-5) {};

  \node [node2] (4el) at (-5,-4) {\text{ $e=4$ }};

  \node (51) at (2.5,-8) {};
  \node (52) at (3.5,-8) {};
  \node (53) at (2.5,-9) {};
  \node (54) at (3.5,-9) {};

  \node [node2] (5el) at (-5,-8) {\text{ $e=5$ }};

  \node (61) at (2.5,-12) {};
  \node (62) at (3.5,-12) {};
  \node (63) at (2.5,-13) {};
  \node (64) at (3.5,-13) {};
  
  \node [node2] (6el) at (-5,-12) {\text{ $e=6$ }};

  \draw (11) -- (12)
        (11) -- (13)
        (11) -- (14)
        
        (21) -- (22)
        (21) -- (23)
        (23) -- (24)

        (31) -- (32)
        (31) -- (33)
        (33) -- (34)
        (31) -- (34)
        
        (41) -- (42)
        (41) -- (43)
        (43) -- (44)
        (42) -- (44)

        (51) -- (52)
        (51) -- (53)
        (53) -- (54)
        (52) -- (54)
        (51) -- (54)
                
        (61) -- (62)
        (61) -- (63)
        (63) -- (64)
        (62) -- (64)
        (61) -- (64)
        (62) -- (63);

    \draw[very thick]
                      (0,-1.2) -- (0.05,-3.8)
                      (5.5,-1.2) -- (0.05,-3.8)
                      (5.5,-1.2) -- (5.5,-3.8)
                      (0.05,-5.2) -- (3,-7.8);
                      
\draw[line width=0.75mm, draw=green,fill=green!50,opacity=0.3] (-1.2,-1.2) rectangle (1.2,0.3);
\draw[very thick, draw=black!50,fill=black!50,opacity=0.08] (4.7,-1.2) rectangle (6.3,0.2);
\draw[very thick, draw=black!50,fill=black!50,opacity=0.08] (-0.7,-5.2) rectangle (0.8,-3.8);
\draw[very thick, draw=green,fill=green!50,opacity=0.3] (4.7,-5.2) rectangle (6.3,-3.8);
\draw[very thick, draw=green,fill=green!50,opacity=0.3] (2.3,-9.2) rectangle (3.7,-7.8);
\draw[very thick, draw=green,fill=green!50,opacity=0.3] (2.3,-13.2) rectangle (3.7,-11.8);;

\draw[very thick, draw=green]
                      (5.5,-5.2) -- (3,-7.8)
                      (3,-9.2) -- (3,-11.8);
\end{tikzpicture}

\caption{The GRAPH of graphs for $n=4$, the optimal graph for a given $e$ is highlighted by green, EDGEs between optimal graphs are colored green.}
\label{fig:5}
\end{center}
\end{figure}

%% file: n4Optimal.tex
\begin{figure}[htbp!]
	\centering
	\pgfplotsset{scaled y ticks=false}
	\pgfplotsset{every axis plot/.append style={line width=1.5pt}}
	\begin{tikzpicture}[scale=0.9]
	\begin{scope}
\begin{axis}[
    width=7cm,
    axis lines = left,
    title = {Errors of optimal graphs for $n=4$},
    xlabel = Number of comparisons ($e$),
    ylabel = $d_{euc}$,
    xmin=3, xmax=6,
    ymin=0, ymax=0.17,
    yticklabel style={
        /pgf/number format/fixed,
        /pgf/number format/precision=2
},
scaled y ticks=false,
    legend style ={ at={(1.03,1)}, 
        anchor=north west, 
        fill=white,align=left},
]

\addplot [
  color=green,
    mark=triangle*
    ]
    coordinates {
    (3,0.1620)(4,0.0970)(5,0.0608)(6,0)
    };

\addplot [
  color=red,
    mark=otimes*
    ]
    coordinates {
    (3,0.1293)(4,0.0771)(5,0.0482)(6,0)
    };

\addplot [
    color=blue,
    mark=square*
    ]
    coordinates {
    (3,0.0918)(4,0.0543)(5,0.0339)(6,0)
    };

\end{axis}
\end{scope}

\begin{scope}[xshift=7.5cm]
\begin{axis}[
    width=7cm,
    axis lines = left,
    title = {Kendall's $\tau$ of optimal graphs for $n=4$},
    xlabel = Number of comparisons ($e$),
    ylabel = Kendall's $\tau$,
    xmin=3, xmax=6,
    ymin=0.5, ymax=1,
    yticklabel style={
        /pgf/number format/fixed,
        /pgf/number format/precision=2
},
scaled y ticks=false,
    legend style ={ at={(1.03,1)}, 
        anchor=north west, 
        fill=white,align=left},
]

\addplot [
  color=green,
    mark=triangle*
    ]
    coordinates {
    (3,0.6164)(4,0.7328)(5,0.8339)(6,1)
    };
\addlegendentry{Strong}

\addplot [
  color=red,
    mark=otimes*
    ]
    coordinates {
    (3,0.6639)(4,0.7705)(5,0.8583)(6,1)
    };
\addlegendentry{Modest}

\addplot [
    color=blue,
    mark=square*
    ]
    coordinates {
    (3,0.7306)(4,0.8216)(5,0.8909)(6,1)
    };
\addlegendentry{Weak}

\end{axis}
\end{scope}
\end{tikzpicture}

\caption{The relation between the number of comparisons ($e$), the errors (Euclidean distances) and Kendall's $\tau$ measures of optimal graphs for $n=4$.}
\label{fig:6}
\end{figure}

%% file: n5graph.tex
\begin{figure}[htbp!]
\centering
\begin{tikzpicture}[every node/.style={circle,inner sep=2pt,draw=black,fill=black!20}, scale=0.8]

\tikzstyle{node2} = [draw=white,fill=white];

 \node (1A) at (0,0.4) {};
  \node (1B) at (-0.5,-1.6) {};
  \node (1C) at (-1,-0.6) {};
  \node (1D) at (0.5,-1.6) {};
  \node (1E) at (1,-0.6) {};
  
  \node (2A) at (4,0.4) {};
  \node (2B) at (3,-0.6) {};
  \node (2C) at (4.5,-1.6) {};
  \node (2D) at (5,-0.6) {};
  \node (2E) at (3.5,-1.6) {};
  
  \node (3A) at (8.5,-1.6) {};
  \node (3B) at (7.5,-1.6) {};
  \node (3C) at (7,-0.6) {};
  \node (3D) at (8,0.4) {};
  \node (3E) at (9,-0.6) {};
  \node [node2] (4el) at (-5,0.4) {$e=4$};
  
   \node (4A) at (-3,-3.3) {};
  \node (4B) at (-3.5,-5.3) {};
  \node (4C) at (-4,-4.3) {};
  \node (4D) at (-2.5,-5.3) {};
  \node (4E) at (-2,-4.3) {};
  
  \node (5A) at (0.5,-3.3) {};
  \node (5B) at (1.5,-4.3) {};
  \node (5C) at (-0.5,-4.3) {};
  \node (5D) at (1,-5.3) {};
  \node (5E) at (0,-5.3) {};
  
   \node (6A) at (4,-3.3) {};
  \node (6B) at (5,-4.3) {};
  \node (6C) at (4.5,-5.3) {};
  \node (6D) at (3,-4.3) {};
  \node (6E) at (3.5,-5.3) {};
  
  \node (7A) at (6.5,-4.3) {};
  \node (7B) at (7,-5.3) {};
  \node (7C) at (7.5,-3.3) {};
  \node (7D) at (8.5,-4.3) {};
  \node (7E) at (8,-5.3) {};
  
  \node (8A) at (10.5,-3.3) {};
  \node (8B) at (11,-5.3) {};
  \node (8C) at (9.5,-4.3) {};
  \node (8D) at (11.5,-4.3) {};
  \node (8E) at (10,-5.3) {};
  \node [node2] (5el) at (-5,-3.3) {$e=5$};
  
  \node (9A) at (-3,-7.2) {};
  \node (9B) at (-3.5,-9.2) {};
  \node (9C) at (-4,-8.2) {};
  \node (9D) at (-2.5,-9.2) {};
  \node (9E) at (-2,-8.2) {};
  
  \node (10A) at (0.5,-7.2) {};
  \node (10B) at (1,-9.2) {};
  \node (10C) at (-0.5,-8.2) {};
  \node (10D) at (1.5,-8.2) {};
  \node (10E) at (0,-9.2) {};
  
   \node (11A) at (3,-8.2) {};
  \node (11B) at (3.5,-9.2) {};
  \node (11C) at (4,-7.2) {};
  \node (11D) at (5,-8.2) {};
  \node (11E) at (4.5,-9.2) {};
  
  \node (12A) at (8.5,-8.2) {};
  \node (12B) at (7,-9.2) {};
  \node (12C) at (7.5,-7.2) {};
  \node (12D) at (8,-9.2) {};
  \node (12E) at (6.5,-8.2) {};
 
  \node (13A) at (10.5,-7.2) {};
  \node (13B) at (11.5,-8.2) {};
  \node (13C) at (11,-9.2) {};
  \node (13D) at (9.5,-8.2) {};
  \node (13E) at (10,-9.2) {};
  \node [node2] (6el) at (-5,-7.2) {$e=6$};
  
    \node (14A) at (-1,-10.7) {};
  \node (14B) at (-2,-11.7) {};
  \node (14C) at (-0.5,-12.7) {};
  \node (14D) at (0,-11.7) {};
  \node (14E) at (-1.5,-12.7) {};
  
   \node (15A) at (2,-10.7) {};
  \node (15B) at (1.5,-12.7) {};
  \node (15C) at (1,-11.7) {};
  \node (15D) at (2.5,-12.7) {};
  \node (15E) at (3,-11.7) {};

  \node (16A) at (5,-12.7) {};
  \node (16B) at (4.5,-11.7) {};
  \node (16C) at (5.5,-10.7) {};
  \node (16D) at (6,-12.7) {};
  \node (16E) at (6.5,-11.7) {};
  
  \node (17A) at (7.5,-11.7) {};
  \node (17B) at (9,-12.7) {};
  \node (17C) at (8.5,-10.7) {};
  \node (17D) at (8,-12.7) {};
  \node (17E) at (9.5,-11.7) {};
  \node [node2] (7el) at (-5,-10.7) {$e=7$};
  
  \node (18A) at (3,-15.2) {};
  \node (18B) at (2.5,-16.2) {};
  \node (18C) at (1.5,-16.2) {};
  \node (18D) at (2,-14.2) {};
  \node (18E) at (1,-15.2) {};
  
  \node (19A) at (6.5,-15.2) {};
  \node (19B) at (5,-16.2) {};
  \node (19C) at (5.5,-14.2) {};
  \node (19D) at (6,-16.2) {};
  \node (19E) at (4.5,-15.2) {};
  \node [node2] (8el) at (-5,-14.2) {$e=8$};
  
  \node (20A) at (4,-17.2) {};
  \node (20B) at (3.5,-19.2) {};
  \node (20C) at (3,-18.2) {};
  \node (20D) at (5,-18.2) {};
  \node (20E) at (4.5,-19.2) {};
  \node [node2] (9el) at (-5,-17.2) {$e=9$};
  
  \node (21A) at (4,-20.2) {};
  \node (21B) at (4.5,-22.2) {};
  \node (21C) at (3,-21.2) {};
  \node (21D) at (5,-21.2) {};
  \node (21E) at (3.5,-22.2) {};
  \node [node2] (10el) at (-5,-20.2) {$e=10$};

  \draw (1A) -- (1B)
        (1A) -- (1C)
        (1A) -- (1D)
        (1A) -- (1E)
        
        (2A) -- (2B)
        (2A) -- (2C)
        (2A) -- (2D)
        (2C) -- (2E)
        
        (3A) -- (3B)
        (3B) -- (3C)
        (3C) -- (3D)
        (3D) -- (3E)
        
        (4A) -- (4B)
        (4A) -- (4C)
        (4A) -- (4D)
        (4A) -- (4E)
        (4B) -- (4D)
        
        (5A) -- (5B)
        (5A) -- (5C)
        (5A) -- (5D)
        (5C) -- (5E)
        (5D) -- (5E)
        
        (6A) -- (6B)
        (6A) -- (6C)
        (6A) -- (6D)
        (6C) -- (6E)
        (6C) -- (6D)
        
        (7B) -- (7E)
        (7A) -- (7C)
        (7A) -- (7D)
        (7C) -- (7E)
        (7C) -- (7D)
        
        (8B) -- (8E)
        (8A) -- (8C)
        (8A) -- (8D)
        (8C) -- (8E)
        (8B) -- (8D)
        
        (9A) -- (9B)
        (9A) -- (9C)
        (9A) -- (9D)
        (9A) -- (9E)
        (9B) -- (9D)
        (9B) -- (9C)
        
        (10A) -- (10B)
        (10A) -- (10C)
        (10A) -- (10D)
        (10C) -- (10E)
        (10D) -- (10E)
        (10B) -- (10E)
        
        (11B) -- (11E)
        (11A) -- (11C)
        (11A) -- (11D)
        (11C) -- (11E)
        (11C) -- (11D)
        (11C) -- (11B)
        
        (12B) -- (12E)
        (12A) -- (12C)
        (12A) -- (12D)
        (12C) -- (12E)
        (12C) -- (12D)
        (12B) -- (12D)
        
        (13A) -- (13B)
        (13A) -- (13C)
        (13A) -- (13D)
        (13C) -- (13E)
        (13C) -- (13D)
        (13B) -- (13D)
        
        (14A) -- (14B)
        (14A) -- (14C)
        (14A) -- (14D)
        (14C) -- (14E)
        (14D) -- (14E)
        (14B) -- (14E)
        (14A) -- (14E)
        
        (15A) -- (15E)
        (15A) -- (15C)
        (15A) -- (15D)
        (15D) -- (15E)
        (15B) -- (15D)
        (15A) -- (15B)
        (15C) -- (15B)
                
        (16A) -- (16B)
        (16A) -- (16C)
        (16A) -- (16D)
        (16C) -- (16E)
        (16C) -- (16D)
        (16B) -- (16D)
        (16C) -- (16B)
        
        (17B) -- (17E)
        (17A) -- (17C)
        (17A) -- (17D)
        (17C) -- (17E)
        (17C) -- (17B)
        (17B) -- (17D)
        (17D) -- (17E)
        
        (18B) -- (18E)
        (18A) -- (18C)
        (18A) -- (18D)
        (18C) -- (18E)
        (18C) -- (18B)
        (18B) -- (18D)
        (18D) -- (18E)
        (18C) -- (18D)
        
        (19B) -- (19E)
        (19A) -- (19C)
        (19A) -- (19D)
        (19C) -- (19E)
        (19B) -- (19A)
        (19B) -- (19D)
        (19A) -- (19E)
        (19C) -- (19D)
        
        (20B) -- (20E)
        (20A) -- (20C)
        (20A) -- (20D)
        (20C) -- (20E)
        (20C) -- (20B)
        (20B) -- (20D)
        (20D) -- (20E)
        (20C) -- (20D)
        (20A) -- (20E)
        
        (21B) -- (21E)
        (21A) -- (21C)
        (21A) -- (21D)
        (21C) -- (21E)
        (21B) -- (21A)
        (21B) -- (21D)
        (21A) -- (21E)
        (21C) -- (21D)
        (21D) -- (21E)
        (21B) -- (21C);

    \draw[very thick] (0,-1.75) -- (-3,-3.1)
                      (4,-1.7) -- (-3,-3.1)
                      (4,-1.7) -- (0.5,-3.1)
                      (4,-1.7) -- (4,-3.1)
                      (4,-1.7) -- (7.5,-3.1)
                      (8,-1.7) -- (0.5,-3.1)
                      (8,-1.7) -- (4,-3.1)
                      (8,-1.7) -- (7.5,-3.1)
                      (8,-1.7) -- (10.5,-3.1)
                      (-3,-5.5) -- (-3,-7)
                      (-3,-5.5) -- (4,-7)
                      (0.5,-5.5) -- (-3,-7)
                      (0.5,-5.5) -- (0.5,-7)
                      (0.5,-5.5) -- (7.5,-7)
                      (0.5,-5.5) -- (10.5,-7)
                      (4,-5.5) -- (-3,-7)
                      (4,-5.5) -- (7.5,-7)
                      (4,-5.5) -- (10.5,-7)
                      (7.5,-5.5) -- (4,-7)
                      (7.5,-5.5) -- (7.5,-7)
                      (7.5,-5.5) -- (10.5,-7)
                      (10.5,-5.5) -- (7.5,-7)
                      (-3,-9.4) -- (-1,-10.4)
                      (-3,-9.4) -- (2,-10.4)
                      (-3,-9.4) -- (5.5,-10.4)
                      (0.5,-9.35) -- (-1,-10.4)
                      (4,-9.4) -- (2,-10.4)
                      (7.5,-9.4) -- (2,-10.4)
                      (7.5,-9.4) -- (8.5,-10.4)
                      (10.5,-9.4) -- (2,-10.4)
                      (10.5,-9.4) -- (5.5,-10.4)
                      (10.5,-9.4) -- (8.5,-10.4)
                      (-1,-12.9) -- (2.5,-14)
                      (2,-12.9) -- (2.5,-14)
                      (5.5,-12.9) -- (2.5,-14)
                      (8.5,-12.9) -- (2.5,-14)
                      (2,-12.9) -- (6,-14)
                      (2,-16.4) -- (4,-17);
                      
\draw[line width=0.75mm, draw=green,fill=green!50,opacity=0.3] (-1.2,-1.75) rectangle (1.2,0.9);
\draw[very thick, draw=black!50,fill=black!50,opacity=0.08] (2.8,-1.7) rectangle (5.2,0.8);
\draw[very thick, draw=black!50,fill=black!50,opacity=0.08] (6.8,-1.7) rectangle (9.2,0.9);
\draw[very thick, draw=black!50,fill=black!50,opacity=0.08] (-4.2,-5.5) rectangle (-1.8,-3.1);
\draw[very thick, draw=black!50,fill=black!50,opacity=0.08] (-0.7,-5.5) rectangle (1.7,-3.1);
\draw[very thick, draw=black!50,fill=black!50,opacity=0.08] (2.8,-5.5) rectangle (5.2,-3.1);
\draw[very thick, draw=black!50,fill=black!50,opacity=0.08] (6.3,-5.5) rectangle (8.7,-3.1);
\draw[line width=0.75mm, draw=green,fill=green!50,opacity=0.3] (9.3,-5.5) rectangle (11.7,-3.1);
\draw[very thick, draw=black!50,fill=black!50,opacity=0.08] (-4.2,-9.4) rectangle (-1.8,-7);
\draw[line width=0.75mm, draw=green,fill=green!50,opacity=0.3] (-0.7,-9.35) rectangle (1.7,-7);
\draw[very thick, draw=black!50,fill=black!50,opacity=0.08] (2.8,-9.4) rectangle (5.2,-7);
\draw[very thick, draw=black!50,fill=black!50,opacity=0.08] (6.3,-9.4) rectangle (8.7,-7);
\draw[very thick, draw=black!50,fill=black!50,opacity=0.08] (9.3,-9.4) rectangle (11.7,-7);
\draw[very thick, draw=black!50,fill=black!50,opacity=0.08] (-2.2,-12.9) rectangle (0.2,-10.4);
\draw[very thick, draw=black!50,fill=black!50,opacity=0.08] (0.8,-12.9) rectangle (3.2,-10.4);
\draw[very thick, draw=black!50,fill=black!50,opacity=0.08] (4.3,-12.9) rectangle (6.7,-10.4);
\draw[line width=0.75mm, draw=green,fill=green!50,opacity=0.3] (7.3,-12.9) rectangle (9.7,-10.4);
\draw[very thick, draw=black!50,fill=black!50,opacity=0.08] (0.8,-16.4) rectangle (3.2,-14);
\draw[line width=0.75mm, draw=green,fill=green!50,opacity=0.3] (4.2,-16.4) rectangle (6.7,-14);
\draw[line width=0.75mm, draw=green,fill=green!50,opacity=0.3] (2.8,-19.4) rectangle (5.2,-17);
\draw[line width=0.75mm, draw=green,fill=green!50,opacity=0.3] (2.7,-22.4) rectangle (5.2,-20);

\draw[very thick, draw=green]
                      (8.5,-12.9) -- (6,-14)
                      (5.5,-16.4) -- (4,-17)
                      (4,-19.4) -- (4,-20)
                      (0.5,-9.35) -- (8.5,-10.4);
\end{tikzpicture}
\caption{The GRAPH of graphs for $n=5$, optimal graphs are highlighted by green.}
 \label{fig:7}
\end{figure}

%% file: n5Optimal.tex
\begin{figure}[htbp!]
	\centering
	\pgfplotsset{scaled y ticks=false}
	\pgfplotsset{every axis plot/.append style={line width=1.5pt}}
	\begin{tikzpicture}[scale=0.9]
	\begin{scope}
\begin{axis}[
    width=7cm,
    axis lines = left,
    title = {Errors of optimal graphs for $n=5$},
    xlabel = Number of comparisons ($e$),
    ylabel = $d_{euc}$,
    xmin=4, xmax=10,
    ymin=0, ymax=0.18,
    yticklabel style={
        /pgf/number format/fixed,
        /pgf/number format/precision=2
},
scaled y ticks=false,
    legend style ={ at={(1.03,1)}, 
        anchor=north west, 
        fill=white,align=left},
]

\addplot [
  color=green,
    mark=triangle*
    ]
    coordinates {
    (4,0.1744)(5,0.1258)(6,0.0988)(7,0.0792)(8,0.0582)(9,0.0365)(10,0)
    };

\addplot [
  color=red,
    mark=otimes*
    ]
    coordinates {
    (4,0.1395)(5,0.1002)(6,0.0785)(7,0.0629)(8,0.0461)(9,0.0289)(10,0)
    };

\addplot [
    color=blue,
    mark=square*
    ]
    coordinates {
    (4,0.0988)(5,0.0707)(6,0.0553)(7,0.0443)(8,0.0324)(9,0.0203)(10,0)
    };

\end{axis}
\end{scope}

\begin{scope}[xshift=7.5cm]
\begin{axis}[
    width=7cm,
    axis lines = left,
    title = {Kendall's $\tau$ of optimal graphs for $n=5$},
    xlabel = Number of comparisons ($e$),
    ylabel = Kendall's $\tau$,
    xmin=4, xmax=10,
    ymin=0.5, ymax=1,
    yticklabel style={
        /pgf/number format/fixed,
        /pgf/number format/precision=2
},
scaled y ticks=false,
    legend style ={ at={(1.03,1)}, 
        anchor=north west, 
        fill=white,align=left},
]

\addplot [
  color=green,
    mark=triangle*
    ]
    coordinates {
    (4,0.5732)(5,0.6546)(6,0.7149)(7,0.7653)(8,0.8198)(9,0.8932)(10,1)
    };
\addlegendentry{Strong}

\addplot [
  color=red,
    mark=otimes*
    ]
    coordinates {
    (4,0.6264)(5,0.7025)(6,0.7564)(7,0.8008)(8,0.8480)(9,0.9108)(10,1)
    };
\addlegendentry{Modest}

\addplot [
    color=blue,
    mark=square*
    ]
    coordinates {
    (4,0.7037)(5,0.7685)(6,0.8127)(7,0.8479)(8,0.8846)(9,0.9329)(10,1)
    };
\addlegendentry{Weak}

\end{axis}
\end{scope}
\end{tikzpicture}

\caption{The relation between the number of comparisons ($e$), the errors (Euclidean distances) and Kendall's $\tau$ measures of optimal graphs for $n=5$.}
\label{fig:8}
\end{figure}

%% file: n6graph.tex
\begin{figure}[htbp!]
\centering
\begin{tikzpicture}[every node/.style={circle,inner sep=2pt,draw=black,fill=black}, scale=0.75]

\tikzstyle{block} = [draw=white,fill=white];

\tikzstyle{node2} = [circle,inner sep=2pt,draw=green,fill=green];

\tikzstyle{node3} = [circle,inner sep=2pt,draw=green!50!orange!50!,fill=green!50!orange!50!];

 \node [node2] (51) at (-1.5,0) {};
  \node (52) at (-1,0) {};
  \node (53) at (-0.5,0) {};
  \node (54) at (0,0) {};
  \node (55) at (0.5,0) {};
  \node (56) at (1,0) {};
  
  \node [block] at (-7,0) (5el) {$e=5$};
  
   \node (61) at (-3,-2) {};
  \node (62) at (-2.5,-2) {};
  \node (63) at (-2,-2) {};
  \node (64) at (-1.5,-2) {};
  \node (65) at (-1,-2) {};
  \node (66) at (-0.5,-2) {};
  \node (67) at (0,-2) {};
  \node (68) at (0.5,-2) {};
  \node (69) at (1,-2) {};
  \node (610) at (1.5,-2) {};
  \node (611) at (2,-2) {};
  \node (612) at (2.5,-2) {};
  \node [node2] (613) at (3,-2) {};
  
  \node [block] at (-7,-2) (6el) {$e=6$};
  
  \node (71) at (-4.5,-4) {};
  \node (72) at (-4,-4) {};
  \node (73) at (-3.5,-4) {};
  \node (74) at (-3,-4) {};
  \node (75) at (-2.5,-4) {};
  \node (76) at (-2,-4) {};
  \node (77) at (-1.5,-4) {};
  \node (78) at (-1,-4) {};
  \node (79) at (-0.5,-4) {};
  \node (710) at (0,-4) {};
  \node [node2] (711) at (0.5,-4) {};
  \node (712) at (1,-4) {};
  \node (713) at (1.5,-4) {};
   \node (714) at (2,-4) {};
  \node (715) at (2.5,-4) {};
  \node (716) at (3,-4) {};
  \node (717) at (3.5,-4) {};
  \node (718) at (4,-4) {};
  \node (719) at (4.5,-4) {};
  
  \node [block] at (-7,-4) (7el) {$e=7$};
  
  \node (81) at (-5.5,-6) {};
  \node (82) at (-5,-6) {};
  \node (83) at (-4.5,-6) {};
  \node (84) at (-4,-6) {};
  \node (85) at (-3.5,-6) {};
  \node (86) at (-3,-6) {};
  \node (87) at (-2.5,-6) {};
  \node (88) at (-2,-6) {};
  \node (89) at (-1.5,-6) {};
  \node (810) at (-1,-6) {};
  \node (811) at (-0.5,-6) {};
  \node (812) at (0,-6) {};
  \node (813) at (0.5,-6) {};
  \node (814) at (1,-6) {};
  \node (815) at (1.5,-6) {};
   \node [node2] (816) at (2,-6) {};
  \node (817) at (2.5,-6) {};
  \node (818) at (3,-6) {};
  \node (819) at (3.5,-6) {};
  \node (820) at (4,-6) {};
  \node (821) at (4.5,-6) {};
  \node (822) at (5,-6) {};
  
  \node [block] at (-7,-6) (8el) {$e=8$};
  
    \node (91) at (-5,-8) {};
  \node (92) at (-4.5,-8) {};
  \node (93) at (-4,-8) {};
  \node (94) at (-3.5,-8) {};
  \node (95) at (-3,-8) {};
  \node (96) at (-2.5,-8) {};
  \node (97) at (-2,-8) {};
  \node (98) at (-1.5,-8) {};
  \node (99) at (-1,-8) {};
  \node (910) at (-0.5,-8) {};
  \node (911) at (0,-8) {};
  \node (912) at (0.5,-8) {};
  \node [node2] (913) at (1,-8) {};
  \node (914) at (1.5,-8) {};
   \node (915) at (2,-8) {};
  \node (916) at (2.5,-8) {};
  \node (917) at (3,-8) {};
  \node (918) at (3.5,-8) {};
  \node (919) at (4,-8) {};
  \node (920) at (4.5,-8) {};
  
  \node [block] at (-7,-8) (9el) {$e=9$};
  
  \node (101) at (-3.5,-10) {};
  \node (102) at (-3,-10) {};
  \node (103) at (-2.5,-10) {};
  \node (104) at (-2,-10) {};
  \node [node2] (105) at (-1.5,-10) {};
  \node (106) at (-1,-10) {};
  \node (107) at (-0.5,-10) {};
  \node (108) at (0,-10) {};
  \node (109) at (0.5,-10) {};
  \node (1010) at (1,-10) {};
  \node (1011) at (1.5,-10) {};
   \node (1012) at (2,-10) {};
  \node (1013) at (2.5,-10) {};
  \node (1014) at (3,-10) {};
  
  \node [block] at (-7,-10) (10el) {$e=10$};
  
  \node (111) at (-2,-12) {};
  \node (112) at (-1.5,-12) {};
  \node (113) at (-1,-12) {};
  \node (114) at (-0.5,-12) {};
  \node (115) at (0,-12) {};
  \node (116) at (0.5,-12) {};
  \node (117) at (1,-12) {};
  \node (118) at (1.5,-12) {};
  \node [node2] (119) at (2,-12) {};
  
   \node [block] at (-7,-12) (11el) {$e=11$};

  \node [node3] (121) at (-1,-14) {};
  \node (122) at (-0.5,-14) {};
  \node (123) at (0,-14) {};
  \node [node3] (124) at (0.5,-14) {};
  \node [node3] (125) at (1,-14) {\tiny E};
  
  \node [block] at (-7,-14) (12el) {$e=12$};
  
   \node [node3] (131) at (-0.5,-16) {};
  \node [node3] (132) at (0.5,-16) {\tiny E};
  
  \node [block] at (-7,-16) (13el) {$e=13$};
  
  \node [node2] (141) at (0,-18) {};
  
   \node [block] at (-7,-18) (14el) {$e=14$};
  
  \node [node2] (151) at (0,-20) {};
  
  \node [block] at (-7,-20) (15el) {$e=15$};
  
  \draw [draw=black!30] (51) -- (61)
        (52) -- (61)
        (52) -- (62)
        (52) -- (63)
        (52) -- (65)
        (53) -- (62)
        (53) -- (63)
        (53) -- (64)
        (53) -- (67)
        (53) -- (610)
        (53) -- (611)
        (53) -- (612)
        (54) -- (63)
        (54) -- (66)
        (54) -- (612)
        (55) -- (65)
        (55) -- (66)
        (55) -- (67)
        (55) -- (68)
        (55) -- (69)
        (55) -- (611)
        (56) -- (66)
        (56) -- (67)
        (56) -- (69)
        (56) -- (610)
        (56) -- (611)
        (56) -- (613)
        (61) -- (71)
        (61) -- (74)
        (62) -- (71)
        (62) -- (72)
        (62) -- (75)
        (62) -- (77)
        (62) -- (710)
        (63) -- (71)
        (63) -- (73)
        (63) -- (75)
        (63) -- (76)
        (63) -- (79)
        (63) -- (710)
        (64) -- (72)
        (64) -- (73)
        (64) -- (711)
        (64) -- (716)
        (64) -- (717)
        (65) -- (74)
        (65) -- (75)
        (65) -- (76)
        (65) -- (77)
        (65) -- (78)
        (65) -- (79)
        (66) -- (75)
        (66) -- (76)
        (66) -- (712)
        (66) -- (713)
        (66) -- (717)
        (67) -- (75)
        (67) -- (711)
        (67) -- (714)
        (67) -- (715)
        (67) -- (717)
        (68) -- (76)
        (68) -- (714)
        (69) -- (77)
        (69) -- (78)
        (69) -- (711)
        (69) -- (712)
        (69) -- (713)
        (69) -- (714)
        (69) -- (718)
        (610) -- (77)
        (610) -- (79)
        (610) -- (715)
        (610) -- (717)
        (610) -- (718)
        (610) -- (719)
        (611) -- (75)
        (611) -- (76)
        (611) -- (78)
        (611) -- (79)
        (611) -- (714)
        (611) -- (715)
        (611) -- (716)
        (611) -- (717)
        (611) -- (718)
        (612) -- (79)
        (612) -- (710)
        (612) -- (717)
        (612) -- (719)
        (613) -- (713)
        (613) -- (715)
        (71) -- (81)
(71) -- (83)
(71) -- (84)
(71) -- (85)
(72) -- (81)
(72) -- (82)
(72) -- (86)
(72) -- (810)
(73) -- (81)
(73) -- (88)
(73) -- (89)
(73) -- (811)
(74) -- (83)
(74) -- (84)
(75) -- (83)
(75) -- (86)
(75) -- (88)
(75) -- (89)
(75) -- (810)
(75) -- (812)
(75) -- (813)
(75) -- (815)
(76) -- (83)
(76) -- (87)
(76) -- (88)
(76) -- (811)
(76) -- (812)
(76) -- (813)
(76) -- (814)
(77) -- (84)
(77) -- (86)
(77) -- (813)
(77) -- (818)
(78) -- (84)
(78) -- (87)
(78) -- (89)
(78) -- (813)
(78) -- (814)
(78) -- (819)
(79) -- (84)
(79) -- (88)
(79) -- (814)
(79) -- (815)
(79) -- (818)
(710) -- (85)
(710) -- (88)
(710) -- (810)
(710) -- (818)
(711) -- (86)
(711) -- (89)
(711) -- (821)
(711) -- (822)
(712) -- (86)
(712) -- (87)
(712) -- (812)
(712) -- (816)
(712) -- (817)
(713) -- (813)
(713) -- (816)
(713) -- (821)
(714) -- (812)
(714) -- (813)
(714) -- (814)
(714) -- (821)
(715) -- (813)
(715) -- (815)
(715) -- (820)
(715) -- (821)
(715) -- (822)
(716) -- (810)
(716) -- (811)
(716) -- (814)
(716) -- (822)
(717) -- (86)
(717) -- (88)
(717) -- (810)
(717) -- (814)
(717) -- (815)
(717) -- (817)
(717) -- (820)
(717) -- (821)
(718) -- (813)
(718) -- (814)
(718) -- (817)
(718) -- (818)
(718) -- (819)
(718) -- (822)
(719) -- (818)
(719) -- (820)
(81) -- (91)
(81) -- (92)
(81) -- (93)
(82) -- (91)
(82) -- (94)
(83) -- (92)
(83) -- (93)
(83) -- (95)
(83) -- (96)
(84) -- (92)
(84) -- (96)
(84) -- (914)
(85) -- (93)
(85) -- (914)
(86) -- (92)
(86) -- (94)
(86) -- (97)
(86) -- (911)
(86) -- (915)
(87) -- (92)
(87) -- (98)
(87) -- (910)
(87) -- (916)
(88) -- (92)
(88) -- (93)
(88) -- (99)
(88) -- (910)
(88) -- (911)
(88) -- (912)
(88) -- (917)
(89) -- (92)
(89) -- (911)
(89) -- (918)
(810) -- (93)
(810) -- (94)
(810) -- (99)
(810) -- (911)
(810) -- (915)
(811) -- (93)
(811) -- (910)
(811) -- (918)
(812) -- (95)
(812) -- (97)
(812) -- (98)
(812) -- (99)
(812) -- (910)
(812) -- (911)
(813) -- (96)
(813) -- (97)
(813) -- (98)
(813) -- (911)
(813) -- (917)
(813) -- (918)
(813) -- (919)
(814) -- (96)
(814) -- (99)
(814) -- (910)
(814) -- (911)
(814) -- (916)
(814) -- (917)
(814) -- (919)
(815) -- (96)
(815) -- (911)
(815) -- (912)
(815) -- (915)
(815) -- (919)
(816) -- (97)
(816) -- (98)
(816) -- (913)
(817) -- (97)
(817) -- (99)
(817) -- (915)
(817) -- (916)
(818) -- (914)
(818) -- (915)
(818) -- (917)
(819) -- (914)
(819) -- (916)
(819) -- (918)
(820) -- (915)
(820) -- (917)
(820) -- (920)
(821) -- (97)
(821) -- (911)
(821) -- (919)
(821) -- (920)
(822) -- (915)
(822) -- (918)
(822) -- (919)
(91) -- (101)
(92) -- (101)
(92) -- (102)
(92) -- (104)
(92) -- (107)
(93) -- (101)
(93) -- (103)
(93) -- (104)
(93) -- (107)
(94) -- (101)
(94) -- (106)
(94) -- (108)
(95) -- (102)
(95) -- (103)
(96) -- (102)
(96) -- (104)
(96) -- (107)
(96) -- (109)
(97) -- (102)
(97) -- (105)
(97) -- (106)
(97) -- (1010)
(97) -- (1012)
(98) -- (102)
(98) -- (105)
(98) -- (1012)
(98) -- (1013)
(99) -- (102)
(99) -- (103)
(99) -- (106)
(99) -- (1011)
(99) -- (1014)
(910) -- (103)
(910) -- (104)
(910) -- (1012)
(910) -- (1013)
(910) -- (1014)
(911) -- (102)
(911) -- (104)
(911) -- (106)
(911) -- (1010)
(911) -- (1011)
(911) -- (1012)
(912) -- (104)
(912) -- (1011)
(913) -- (105)
(914) -- (107)
(915) -- (107)
(915) -- (108)
(915) -- (1010)
(915) -- (1011)
(916) -- (107)
(916) -- (1012)
(916) -- (1014)
(917) -- (107)
(917) -- (1010)
(917) -- (1011)
(917) -- (1013)
(918) -- (107)
(918) -- (1012)
(919) -- (109)
(919) -- (1010)
(919) -- (1011)
(919) -- (1012)
(920) -- (1010)
(101) -- (111)
(101) -- (113)
(102) -- (111)
(102) -- (112)
(102) -- (114)
(102) -- (115)
(103) -- (111)
(103) -- (115)
(103) -- (116)
(104) -- (111)
(104) -- (114)
(104) -- (115)
(105) -- (112)
(105) -- (119)
(106) -- (111)
(106) -- (112)
(106) -- (117)
(106) -- (118)
(107) -- (113)
(107) -- (114)
(107) -- (115)
(108) -- (113)
(108) -- (117)
(109) -- (114)
(1010) -- (114)
(1010) -- (117)
(1010) -- (119)
(1011) -- (114)
(1011) -- (115)
(1011) -- (117)
(1012) -- (114)
(1012) -- (115)
(1012) -- (118)
(1012) -- (119)
(1013) -- (115)
(1013) -- (119)
(1014) -- (115)
(1014) -- (116)
(1014) -- (118)
        (111) -- (121)
        (111) -- (122)
        (111) -- (123)
        (112) -- (121)
        (112) -- (124)
        (113) -- (122)
        (114) -- (122)
        (114) -- (124)
        (115) -- (122)
        (115) -- (123)
        (115) -- (124)
        (116) -- (123)
        (117) -- (122)
        (117) -- (124)
        (117) -- (125)
        (118) -- (123)
        (118) -- (124)
        (119) -- (124)
        (121) -- (131)
        (122) -- (131)
        (122) -- (132)
        (123) -- (131)
        (124) -- (131)
        (124) -- (132)
        (125) -- (132)
        (131) -- (141)
        (132) -- (141)
        (141) -- (151);
        
\draw [very thick, draw=green]
        (816) -- (913)
        (913) -- (105)
        (105) -- (119)
        (132) -- (141)
        (141) -- (151);

\draw [very thick, draw=orange]
        (613) -- (713)
        (713) -- (816)
        (119) -- (124)
        (124) -- (132);

\draw [very thick, draw=blue]
        (125) -- (132);

\end{tikzpicture}

\caption{\footnotesize The GRAPH of graphs for $n=6$, optimal graphs (=NODEs) are colored green, EDGEs between optimal graphs are colored green, too. For $e=12$ and $e=13$ some of the Kendall's $\tau$ values are smaller than the margin of error, the graphs (=NODEs) that are at least second according to at least one indicator are denoted by lighter green color, while $E$ denotes the unique optimal graphs (=NODEs) according to the Euclidean distance for a given $e$. These latter graphs (=NODEs) are connected by a blue EDGE, while the path including as many of the optimal graphs (=NODEs) as possible, is complemented by orange EDGEs besides the green ones.}
\label{fig:9}
\end{figure}

%% file: MKE_PTE.tex
\begin{figure}[htbp!]
\centering
\begin{tikzpicture}[every node/.style={circle,inner sep=2pt,draw=black,fill=black},scale=0.9]

\tikzstyle{block} = [circle,inner sep=2pt,draw=white,fill=white];


  \node (651) at (-8,-1.5) {};
  \node (652) at (-7,-2.5) {};
  \node (653) at (-9,-2.5) {};
  \node (654) at (-7,-3.5) {};
  \node (655) at (-9,-3.5) {};
  \node (656) at (-8,-4.5) {};
  
  \node [block] at (-8,0) (5) {\small $n=6, e=5$};
  
  \node (661) at (-2.75,-1.5) {};
  \node (662) at (-1.75,-2.5) {};
  \node (663) at (-3.75,-2.5) {};
  \node (664) at (-1.75,-3.5) {};
  \node (665) at (-3.75,-3.5) {};
  \node (666) at (-2.75,-4.5) {};
  
  \node [block] at (-2.75,0) (6) {\small $n=6, e=6$};
  
  \node (671) at (2.5,-1.5) {};
  \node (672) at (3.5,-2.5) {};
  \node (673) at (1.5,-2.5) {};
  \node (674) at (3.5,-3.5) {};
  \node (675) at (1.5,-3.5) {};
  \node (676) at (2.5,-4.5) {};
  
  \node [block] at (2.5,0) (7) {\small $n=6, e=7$};
  
  \node (681) at (-10,-7.5) {};
  \node (682) at (-9,-8.5) {};
  \node (683) at (-11,-8.5) {};
  \node (684) at (-9,-9.5) {};
  \node (685) at (-11,-9.5) {};
  \node (686) at (-10,-10.5) {};
  
  \node [block] at (-10,-6) (8) {\small $n=6, e=8$};
  
  \node (691) at (-4.75,-7.5) {};
  \node (692) at (-3.75,-8.5) {};
  \node (693) at (-5.75,-8.5) {};
  \node (694) at (-3.75,-9.5) {};
  \node (695) at (-5.75,-9.5) {};
  \node (696) at (-4.75,-10.5) {};
  
  \node [block] at (-4.75,-6) (9) {\small $n=6, e=9$};
  
  \node (6101) at (0,-7.5) {};
  \node (6102) at (-1,-8.5) {};
  \node (6103) at (1,-8.5) {};
  \node (6104) at (-1,-9.5) {};
  \node (6105) at (1,-9.5) {};
  \node (6106) at (0,-10.5) {};
  
  \node [block] at (0,-6) (10) {\small $n=6, e=10$};
  
  \node (6111) at (4.75,-7.5) {};
  \node (6112) at (3.75,-8.5) {};
  \node (6113) at (5.75,-8.5) {};
  \node (6114) at (3.75,-9.5) {};
  \node (6115) at (5.75,-9.5) {};
  \node (6116) at (4.75,-10.5) {};
  
  \node [block] at (4.75,-6) (11) {\small $n=6, e=11$};
  
  \node (6121) at (-10,-13.5) {};
  \node (6122) at (-9,-14.5) {};
  \node (6123) at (-11,-14.5) {};
  \node (6124) at (-9,-15.5) {};
  \node (6125) at (-11,-15.5) {};
  \node (6126) at (-10,-16.5) {};
  
  \node [block] at (-10,-12) (12) {\small $n=6, e=12$};
  
  \node (6131) at (-3.75,-14.5) {};
  \node (6132) at (-5.75,-14.5) {};
  \node (6133) at (-4.75,-13.5) {};
  \node (6134) at (-3.75,-15.5) {};
  \node (6135) at (-5.75,-15.5) {};
  \node (6136) at (-4.75,-16.5) {};
  
  \node [block] at (-4.75,-12) (13) {\small $n=6, e=13$};
  
  \node (6141) at (1,-14.5) {};
  \node (6142) at (-1,-14.5) {};
  \node (6143) at (0,-13.5) {};
  \node (6144) at (1,-15.5) {};
  \node (6145) at (-1,-15.5) {};
  \node (6146) at (0,-16.5) {};
  
  \node [block] at (0,-12) (9) {\small $n=6, e=14$};
  
   \node (6151) at (5.75,-14.5) {};
  \node (6152) at (3.75,-14.5) {};
  \node (6153) at (4.75,-13.5) {};
  \node (6154) at (5.75,-15.5) {};
  \node (6155) at (3.75,-15.5) {};
  \node (6156) at (4.75,-16.5) {};
  
  \node [block] at (4.75,-12) (9) {\small $n=6, e=15$};
  
  \draw [draw=black!50]
                        (651) -- (652)
                        (651) -- (653)
                        (651) -- (654)
                        (651) -- (655)
                        (651) -- (656)
                        (661) -- (662)
                        (661) -- (663)
                        (662) -- (664)
                        (663) -- (665)
                        (664) -- (666)
                        (665) -- (666)
                        (671) -- (672)
                        (671) -- (673)
                        (672) -- (674)
                        (673) -- (675)
                        (674) -- (675)
                        (674) -- (676)
                        (673) -- (676)
                        (681) -- (682)
                        (681) -- (683)
                        (682) -- (684)
                        (683) -- (685)
                        (684) -- (686)
                        (685) -- (686)
                        (682) -- (685)
                        (683) -- (684)
                        (691) -- (692)
                        (691) -- (693)
                        (692) -- (694)
                        (693) -- (695)
                        (694) -- (696)
                        (695) -- (696)
                        (692) -- (695)
                        (693) -- (694)
                        (691) -- (696)
                        (6101) -- (6102)
                        (6101) -- (6103)
                        (6102) -- (6104)
                        (6103) -- (6105)
                        (6104) -- (6106)
                        (6105) -- (6106)
                        (6102) -- (6105)
                        (6103) -- (6104)
                        (6101) -- (6106)
                        (6101) -- (6104)
                        (6111) -- (6112)
                        (6111) -- (6113)
                        (6112) -- (6114)
                        (6113) -- (6115)
                        (6114) -- (6116)
                        (6115) -- (6116)
                        (6112) -- (6115)
                        (6113) -- (6114)
                        (6111) -- (6116)
                        (6111) -- (6114)
                        (6112) -- (6116)
                        (6121) -- (6122)
                        (6121) -- (6123)
                        (6122) -- (6124)
                        (6122) -- (6126)
                        (6123) -- (6125)
                        (6124) -- (6126)
                        (6125) -- (6126)
                        (6121) -- (6124)
                        (6121) -- (6125)
                        (6122) -- (6123)
                        (6124) -- (6125)
                        (6123) -- (6126)
                        (6131) -- (6132)
                        (6131) -- (6133)
                        (6132) -- (6134)
                        (6132) -- (6136)
                        (6133) -- (6135)
                        (6134) -- (6136)
                        (6135) -- (6136)
                        (6132) -- (6135)
                        (6133) -- (6134)
                        (6131) -- (6136)
                        (6131) -- (6134)
                        (6132) -- (6133)
                        (6134) -- (6135)
                        (6141) -- (6142)
                        (6141) -- (6143)
                        (6142) -- (6144)
                        (6143) -- (6145)
                        (6144) -- (6146)
                        (6145) -- (6146)
                        (6142) -- (6145)
                        (6143) -- (6144)
                        (6141) -- (6146)
                        (6141) -- (6144)
                        (6142) -- (6146)
                        (6142) -- (6143)
                        (6143) -- (6146)
                        (6144) -- (6145)
                        (6151) -- (6152)
                        (6151) -- (6153)
                        (6152) -- (6154)
                        (6153) -- (6155)
                        (6154) -- (6156)
                        (6155) -- (6156)
                        (6152) -- (6155)
                        (6153) -- (6154)
                        (6151) -- (6156)
                        (6151) -- (6154)
                        (6151) -- (6155)
                        (6152) -- (6156)
                        (6152) -- (6153)
                        (6153) -- (6156)
                        (6154) -- (6155);
\draw [very thick, draw=green]
                        (691) -- (696)
                        (6101) -- (6104)
                        (6112) -- (6116)
                        (6143) -- (6146)
                        (6151) -- (6155);
                        
\draw [very thick, draw=blue]
                        (6132) -- (6134);

\end{tikzpicture}

\caption{The optimal graphs related to the green NODEs in Figure \ref{fig:9}. The second row shows a partial optimal filling sequence corresponding to the one in Figure \ref{fig:9}, these graphs can be reached from each other. The additional comparisons are highlighted in every step according to the color of the corresponding EDGE in Figure \ref{fig:9}.}
\label{fig:10}
\end{figure}

%% file: n6Optimal.tex
\begin{figure}[ht!]
	\centering
	\pgfplotsset{scaled y ticks=false}
	\pgfplotsset{every axis plot/.append style={line width=1.5pt}}
	\begin{tikzpicture}[scale=0.9]
	\begin{scope}
\begin{axis}[
    width=7cm,
    axis lines = left,
    title = {Errors of optimal graphs for $n=6$},
    xlabel = Number of comparisons ($e$),
    ylabel = $d_{euc}$,
    xmin=5, xmax=15,
    ymin=0, ymax=0.18,
    yticklabel style={
        /pgf/number format/fixed,
        /pgf/number format/precision=2
},
scaled y ticks=false,
    legend style ={ at={(1.03,1)}, 
        anchor=north west, 
        fill=white,align=left},
]

\addplot [
  color=green,
    mark=triangle*
    ]
    coordinates {
    (5,0.1773)(6,0.1418)(7,0.1179)(8,0.1012)(9,0.0823)(10,0.0729)(11,0.0625)(12,0.0500)(13,0.0390)(14,0.0245)(15,0)
    };

\addplot [
  color=red,
    mark=otimes*
    ]
    coordinates {
    (5,0.1419)(6,0.1129)(7,0.0938)(8,0.0803)(9,0.0654)(10,0.0578)(11,0.0495)(12,0.0396)(13,0.0308)(14,0.0194)(15,0)
    };

\addplot [
    color=blue,
    mark=square*
    ]
    coordinates {
    (5,0.1007)(6,0.0798)(7,0.0662)(8,0.0566)(9,0.0461)(10,0.0407)(11,0.0349)(12,0.0278)(13,0.0216)(14,0.0135)(15,0)
    };

\end{axis}
\end{scope}

\begin{scope}[xshift=7.5cm]
\begin{axis}[
    width=7cm,
    axis lines = left,
    title = {Kendall's $\tau$ of optimal graphs for $n=6$},
    xlabel = Number of comparisons ($e$),
    ylabel = Kendall's $\tau$,
    xmin=5, xmax=15,
    ymin=0.5, ymax=1,
    yticklabel style={
        /pgf/number format/fixed,
        /pgf/number format/precision=2
},
scaled y ticks=false,
    legend style ={ at={(1.03,1)}, 
        anchor=north west, 
        fill=white,align=left},
]

\addplot [
  color=green,
    mark=triangle*
    ]
    coordinates {
    (5,0.5464)(6,0.6035)(7,0.6550)(8,0.6951)(9,0.7399)(10,0.7679)(11,0.7986)(12,0.8382)(13,0.8758)(14,0.9259)(15,1)
    };
\addlegendentry{Strong}

\addplot [
  color=red,
    mark=otimes*
    ]
    coordinates {
    (5,0.6049)(6,0.6586)(7,0.7051)(8,0.7408)(9,0.7800)(10,0.8046)(11,0.8310)(12,0.8642)(13,0.8966)(14,0.9386)(15,1)
    };
\addlegendentry{Modest}

\addplot [
    color=blue,
    mark=square*
    ]
    coordinates {
    (5,0.6870)(6,0.7338)(7,0.7722)(8,0.8009)(9,0.8325)(10,0.8517)(11,0.8725)(12,0.8977)(13,0.9226)(14,0.9547)(15,1)
    };
\addlegendentry{Weak}

\end{axis}
\end{scope}
\end{tikzpicture}

\caption{The relation between the number of comparisons ($e$), the errors (Euclidean distances) and Kendall's $\tau$ measures of optimal graphs for $n=6$.}
\label{fig:11}
\end{figure}

%% file: CR_lognorm.tex
\begin{figure}[ht!]
\renewcommand{\thefigure}{A\arabic{appendfigure}}
	\centering
	\pgfplotsset{scaled y ticks=false}
	\pgfplotsset{every axis plot/.append style={line width=1.5pt}}
	\begin{tikzpicture}
	\begin{scope}
\begin{axis}[box plot width=1.5mm,
title = {Weak perturbation},
ylabel = Consistency Ratio (CR),
ymin=0,ymax=0.5,
xtick={5,6,7,8,9,10},
width=5cm,
height=8cm]
\boxplot[forget plot, blue,box plot whisker bottom index=1,
    box plot whisker top index=5,
    box plot box bottom index=2,
    box plot box top index=4,
    box plot median index=3]{weaklog.dat};

\end{axis}
\end{scope}

\begin{scope}[xshift=5.5cm]
\begin{axis}[box plot width=1.5mm,
ymin=0,ymax=0.5,
xtick={5,6,7,8,9,10},
title = {Modest perturbation},
xlabel = Number of alternatives/criteria ($n$),
width=5cm,
height=8cm]

\boxplot[forget plot, red,box plot whisker bottom index=1,
    box plot whisker top index=5,
    box plot box bottom index=2,
    box plot box top index=4,
    box plot median index=3]{modestlog.dat};

\end{axis}
\end{scope}

\begin{scope}[xshift=11cm]
\begin{axis}[box plot width=1.5mm,
ymin=0,ymax=0.5,
xtick={5,6,7,8,9,10},
title = {Strong perturbation},
width=5cm,
height=8cm]

\boxplot[forget plot, green,box plot whisker bottom index=1,
    box plot whisker top index=5,
    box plot box bottom index=2,
    box plot box top index=4,
    box plot median index=3]{stronglog.dat};

\end{axis}
\end{scope}

\end{tikzpicture}

\caption{The relation between CR and the element-wise perturbation applying lognormal errors via Box plots. Each Box plot is based on 1000 randomly generated perturbed PCMs, and shows the minimum, maximum, and the first, second (median), and third quartile of the gained CR values. \color{black}}
\label{fig:append1}
\end{figure}